%% file: ARAS-theo.tex
\documentclass[10pt,a4paper]{siamltex}
\usepackage[english]{babel}

\usepackage{latexsym}
\usepackage{amsmath}
\usepackage{amssymb}
\usepackage{amsfonts}
\usepackage{algorithm,algorithmic}
\usepackage{stmaryrd} 

\usepackage{enumerate}

\usepackage{array}
\usepackage{tikz}
\usepackage{pgf}
\usetikzlibrary{patterns}
\usetikzlibrary{arrows}
\usetikzlibrary{snakes}

\usepackage{times}
\newtheorem{remark}{\it Remark\/}
\begin{document}

\title{ARAS: fully algebraic two-level domain decomposition preconditioning technique with approximation on coarse interfaces}

\author{Thomas Dufaud$^{1}$ and Damien Tromeur-Dervout$^{1}$ \\ $^1$ Universit\'e de Lyon, universit\'e Lyon 1, CNRS, institut Camille-Jordan, \\43, boulevard du 11 Novembre 1918, 69622 Villeurbanne, France}

\maketitle 
\begin{abstract}
This paper focuses on the development of a two-level preconditioner based on a fully algebraical enhancement of a Schwarz domain decomposition method. We consider the purely divergence of a Restricted Additive Scwharz iterative process leading to the preconditioner developped by X.-C. Cai and M. Sarkis in SIAM Journal of Scientific Computing, Vol. 21 no. 2, 1999. The convergence of vectorial sequence of traces of this process on the artificial interface can be accelerated by an Aitken acceleration technique as proposed in the work of M. Garbey and D. Tromeur-Dervout,  in International Journal for Numerical Methods in Fluids, Vol. 40, no. 12,2002. We propose a formulation of the Aitken-Schwarz technique as a preconditioning technique called Aitken-RAS \footnote{This paper extends the proposition of the ARAS preconditioning technique published in T. Dufaud and D. Tromeur-Dervout, Aitken's acceleration of the Resctricted Additive Schwarz preconditioning using coarse approximations on the interface. C. R. Math. Acad. Sci. Paris, Vol. 348, no. 13-14, pages 821-824, 2010, by developing the building of the preconditioner and its theoretical properties. Moreover, it focuses on a fully algebraic technique based on SVD to approximate the solutions. Finally, results on industrial linear systems are provided.}. The Aitken acceleration is performed in a reduced space to save computing or permit fully algebraic computation of the accelerated solution without knowledge of the underlying equations. A convergence study of the Aitken-RAS preconditioner is proposed also application on industrial problem.
\end{abstract} 

\section{Introduction}
\label{section-introduction}

The convergence rate of a Krylov method such as GCR  and GMRES, developed by Einsenstat \& al \cite{MR694523}, to solve a linear system 
$Au=f, \, A=(a_{ij}) \in \mathbb{R}^{m \times m}, u \in \mathbb{R}^m, b \in  \mathbb{R}^m$, depends on
the matrix eigenvalues distribution. They derived the convergence rate of the GCR and GMRES methods as for the GCR:
\begin{eqnarray}
||r_{i}||_2 \leq [1-\frac{\lambda_{min}(M)^2}{\lambda_{min}(M)\lambda_{max}(M) + \rho(R)^2}]^{i/2} ||r_0||_2 \leq  [1-\frac{1}{\kappa(M)}]^{i/2} ||r_0||_2 \label{convergenceGCR}
\end{eqnarray}
where $M=(A+A^t)/2$ and $R=(A-A^t)/2$. The GMRES method is mathematically equivalent to the ORTHORES algorithm developed by Young and Cea \cite{MR591431}. 
Its convergence rate  follows the same formula as \eqref{convergenceGCR} if $A$ is positive real. Otherwise when $A$ is diagonalisable
$A=X \Lambda X^{-1}$ , the convergence of the GMRES method depends on the distribution of the eigenvalues and the condition number $\kappa(X)$  as follows.
Let $\lambda_1,\ldots,\lambda_\mu$ be the eigenvalues of $A$ with a  non-positive real part, and let $\Lambda_{\mu 1},\ldots,\lambda_n$  those with 
a positive real part belonging to the circle centred in $C>0$ with a radius $R$ with $C>R$.  Then the GMRES convergence rate can be written as:
\begin{eqnarray}
||r_{i 1}||_2 \leq \kappa(X) \left[ \frac{D}{d}\right]^\mu \left[ \frac{R}{C}\right]^{i-\mu}   ||r_0||_2 \label{convergenceGMRES}
\end{eqnarray}
with $D=\max_{i=1,\mu; j=\mu 1,n} |\lambda_i-\lambda_j|$ and $d=\min_{i=1,\mu}|\lambda_i|$.\eqref{convergenceGCR} and \eqref{convergenceGMRES},
and, generally speaking, decreases when the condition  number $\kappa_2(A)=||A||_2 ||A^{-1}||_2$ 
of the non-singular matrix $A$ increases. This implies the need to reduce the scattering of the eigenvalues distribution 
in the complex plane in order to improve the convergence rate. This is the goal of a preconditioning technique.
The left-preconditioning techniques consist to solve $M^{-1} A u = M^{-1} f$ such that $\kappa_2(M^{-1}A) << \kappa_2(A)$. In this work we focus on the Schwarz preconditioning techniques and the preconditioning techniques that are related to the  Schur complement of the matrix $A$. 

Let us first recall some state of the art about the Schwarz and Aitken-Schwarz solvers, and preconditioners based on the Restricted Additive Schwarz (RAS).

\subsection{State of the art of Schwarz and Aitken techniques}
First, lets us consider the Generalized Schwarz Alternating Method introduced by Engquist and Zao \cite{MR1644668} 
that gathers several Schwarz techniques (see Quarteroni and Valli \cite{bookQuarteroni}). For sake of simplicity, let us consider the case where the whole domain $\Omega$ is decomposed into two sub-domains  $\Omega_1$ and $\Omega_2$, with overlapping or not, defining  two  artificial boundaries $\Gamma_1$, $\Gamma_2$. Let $\Omega_{11}=\Omega_1 \backslash\Omega_2$, $\Omega_{22}=\Omega_2 \backslash \Omega_1$ if there is an overlap. Let $L(x)$ be the continuous operator associated with the discrete operator $A$.
It  can be written in its multiplicative version as:
\begin{algorithm}[H]
\caption{GSAM: Multiplicative version}
\label{algorithm-GSAM-multiplicative}
\begin{algorithmic}[1]
\STATE DO until convergence
\STATE Solve
\begin{eqnarray}
  L(x)u_1^{2n+1}(x) &=& f(x),\; \forall x \in \Omega_1, \\
  u_1^{2n+1}(x) &=& g(x),\;\forall x \in \partial \Omega_1 \backslash \Gamma_1, \\
  \Lambda_1 u_1^{2n+1} &+& \lambda_1 \frac{\partial u_1^{2n+1}(x)}{\partial
  n_1} =\Lambda_1 u_2^{2n} + \lambda_1 \frac{\partial u_2^{2n}(x)}{\partial
  n_1}, \; \forall x \in \Gamma_1
\end{eqnarray}
\STATE Solve
\begin{eqnarray}
  L(x)u_2^{2n+2}(x) &=& f(x),\; \forall x \in \Omega_2, \\
  u_2^{2n+2}(x) &=& g(x),\;\forall x \in \partial \Omega_2 \backslash \Gamma_2, \\
  \Lambda_2 u_2^{2n+2} &+& \lambda_2 \frac{\partial u_2^{2n+2}(x)}{\partial
  n_2} =\Lambda_2 u_1^{2n+1} + \lambda_2 \frac{\partial u_1^{2n+1}(x)}{\partial
  n_2}, \; \forall x \in \Gamma_2.
\end{eqnarray}
\STATE Enddo
\end{algorithmic}
\end{algorithm}
where $\Lambda_i$ are some operators and $\lambda_i$ are constants.

According to the specific choice of the operators $\Lambda_i$ and the values of scalars $\lambda_i$,
we obtain the family of Schwarz domain decomposition techniques:
\begin{table}[H]
\centering
\begin{tabular}{|c|c|c|c|c|c|}
\hline
Overlap & $\Lambda_1$ &$\Lambda_2$ &$\lambda_1$ &$\lambda_2$ & Method \\
\hline
yes & $Id$ & $Id$ & $0$ & $0$ & Schwarz \\
\hline
yes &  $Id$ & $Id$ & $\alpha$ & $\alpha$ & ORAS (St Cyr \& al \cite{MR2357620}) \\
\hline
No & $Id$ & $0$ & $0$ & $1$ &  Neumann-Dirichlet (Marini \& Quarteronni \cite{MR998911}) \\
\hline
No &  $Id$ & $Id$ & $1$ & $1$ & Modified Schwarz (Lions \cite{MR972510}) \\
\hline
\end{tabular}
\caption{Derived methods obtained from the specific choices of the operators $\Lambda_i$ and the values of scalars $\lambda_i$ in the GSAM.}
\label{table-GSAM-derive}
\end{table}
If $\Lambda_1=\Lambda_2=I$ and $\lambda_1=\lambda_2=0$
then the above multiplicative version is the classical
Multiplicative Schwarz. If $\Lambda_1=\Lambda_2=constant$ and $
\lambda_1=\lambda_2=1$ then it is the modified Schwarz proposed by
Lions in \cite{MR1064345}.

Engquist and Zao \cite{MR1644668}  showed that with an appropriate choice of the operators $\Lambda_i$
this domain decomposition method  converges in two iterations.  They established the proposition that follows:
\begin{proposition}If $\Lambda_1$ (or $\Lambda_2$) is the Dirichlet to Neumann
operator at the artificial boundary $\Gamma_1$ (or $\Gamma_2$) for
the corresponding homogeneous PDE in $\Omega_2$ (or $\Omega_1$)
with homogeneous boundary condition on $\partial \Omega_2 \cap
\partial \Omega$ (or  $\partial \Omega_1 \cap
\partial \Omega$) then the Generalized Schwarz Alternating method
converges in two steps.
\end{proposition}

The GSAM method converges in two steps if the Dirichlet-Neumann operators $\Lambda_i$, $i=1,2$, are available. 
These operators are not local  to a sub-domain but  they link up together all the sub-domains. In practice, 
some approximations defined algebraically of these operators are used (see Chevalier \& Nataf \cite{MR1641727}, 
Gander \& al \cite{MR1924414}, Gerardo-Giorda \& Nataf \cite{MR2189549}). 

In the Aitken-Schwarz methodology introduced by Garbey \& Tromeur-Dervout \cite{mg_dtd_dd12,mg_dtd_fd}, only the convergence property 
of the Schwarz  method is used. Consequently, no direct approximation of the Dirichlet-Neumann map is used, 
but an approximation of the operator of error linked to this Dirichlet-Neumann map is performed.
This Aitken-Schwarz methodology is based on the purely linear convergence for 
the Schwarz Alternating method when the local operators are linear operators.
\begin{definition}\label{def-aitken}
Let $\left(\left(u^{k}_{i}\right)_{i=1,...,n} = u^k \right)_{k \in \mathbb{N}}$ be a vectorial sequence converging toward  $\left(\xi\right)_{i=1,...,n}=\xi$ purely linearly if
\begin{equation}
u^{k+1} - \xi = P \left( u^k - \xi \right) 
\end{equation}
where $P \in \mathbb{R}^{n \times n}$ is a constant error's transfer operator independent of $k$ and non-singular.
\end{definition}
We assume that there exists a norm $||.||$ such as $||P||<1$. Then $ P $ and $ \xi $ can be determined from $ n+1 $ iterations , using the equations:
\begin{equation}
\left(u^{k+1}-u^{k}, ..., u^2-u^1\right) = P\left( u^{k}-u^{k-1}, ..., u^1-u^0 \right)
\end{equation}
So, if $\left( u^{n}-u^{n-1}, ..., u^1-u^0 \right)$ is non-singular $ P$ can be written as :
\begin{equation}
P = \left(u^{n+1}-u^{n}, ..., u^2-u^1\right)\left( u^{n}-u^{n-1}, ..., u^1-u^0 \right)^{-1}
\end{equation}

Then if $ ||P||<1 $, $(I - P)$ is non singular and it is possible to compute $\xi$ as
\begin{equation}\label{eq:aitken-formula}
\xi = \left( I - P \right)^{-1}\left(u^{n+1}-Pu^{n} \right)
\end{equation}
For domain decomposition methods, the vectorial sequences $u^n$ corresponds to the iterated solution at the subdomains artificial interfaces. 
To apply directly the Aitken's acceleration in the vectorial case, we have to construct the matrix $P$ or an approximation of it, and to apply the Aitken's acceleration \eqref{eq:aitken-formula}. Algorithm \ref{algorithm-physical-space} describes the acceleration written in the canonical base of $\mathbb{R}^n$ ("physical space").

\begin{algorithm}[H]
\caption{Vectorial Aitken's acceleration in the physical space}
\label{algorithm-physical-space}
\begin{algorithmic}[1]
\REQUIRE {${\cal{G}}:   \mathbb{R}^n \rightarrow \mathbb{R}^n$ an iterative method having a pure linear convergence}
\REQUIRE {$(u^k)_{1 \leq k \leq n+1}$, \; $n+1$ successive iterates of ${\cal{G}}$ starting from an arbitrary initial guess $u^0$}
\STATE Form $E^k=u^{k+1}-u^k, \; {0 \leq k \leq n}$
\IF {$\left[ E^{n-1},\ldots,E^0 \right]$ is invertible} 
\STATE $P=\left[ E^{n},\ldots,E^1 \right] \left[ E^{n-1},\ldots,E^0 \right]^{-1}$
\STATE  $u^\infty = (I_n-P)^{-1}(u^{n+1}-P u^{n})$
\ENDIF
\end{algorithmic}
\end{algorithm}
The drawback of Algorithm \ref{algorithm-physical-space} is to be limited to a sequence of small vector size because it needs a number of iterations related to the vector size $n$. In order to overcome this difficulty, some approximation of the error transfer operator $P$ is proposed using some coarse approximation spaces to represent the solution. Garbey \cite{Garbey1} proposed to write the solution in the eigenbasis associated to the part of the separable operator associated to the direction parallel to the artificial interfaces, or with a coarse approximation of the  sinus or cosinus expansion of the solution at the interface. Tromeur-Dervout \cite{pareng09SVD} proposed to build an approximation space based on the Singular Value Decomposition (SVD) of the interface solutions of the Schwarz iterates. This last approximation override any constraints about separability of the linear operator $A$  and mesh considerations. 
 
Nevertheless this last techniques fail when: 
\begin{itemize}
\item the iterative process based on domain decomposition diverges too fast,
\item the local solutions are inaccurately solved leading to a less numerically efficient acceleration by the Aitken's process.
\end{itemize}

In such cases the Aitken-Schwarz method as solver is no longer suitable should be considered as preconditioner of a Krylov method. The purpose of this paper is to detail and to extend the preconditioners based on Schwarz domain decomposition accelerated by Aitken's techniques developped by Dufaud \& Tromeur-Dervout \cite{ARASCras} with building and approximation of the matrix $P$  arising from the SVD approximation of the Schwarz interface solutions.  

\subsection{State of the art of preconditioners based on RAS}
The techniques  of preconditioning that are based on  domain decomposition of Schwarz's type have been widely developed  this last decade and accelerated multiplicative Schwarz has been 
"{\em a consistently good performer}"  as said Cai \& al \cite{MR1307330}. The first type of  domain decomposition preconditioning to appear was  domain 
decomposition based on substructuring technique of Bramble \& al \cite{MR842125} followed by the Additive Schwarz (AS) preconditioning of 
Dryja \& Widlund \cite{MR1038100}, Gropp \& Keyes \cite{MR1145180}. It is built from the adjacency graph $G=(W,E)$ of $A$, where $W=\lbrace 1,2,...,m \rbrace$ and  $E=\lbrace (i,j) : a_{ij} \neq 0 \rbrace$  are the edges and vertices of $G$.  Starting with a non-overlapping partition $W=\cup_{i=1}^p W_{i,0}$ and $\delta \ge 0$ given, 
the overlapping  partition $\lbrace W_{i,\delta} \rbrace$ is  obtained defining $p$ partitions $ W_{i,\delta} \supset W_{i,\delta-1}$ by including all the immediate neighbouring vertices of the vertices in the partition $W_{i,\delta-1}$. Then the  restriction operator $R_{i,\delta}: W \rightarrow W_{i,\delta}$
defines the  local operator $A_{i,\delta} =R_{i,\delta} A R_{i,\delta}^{T} , A_{i,\delta} \in \mathbb{R}^{m_{i,\delta}\times m_{i,\delta}} $  on $W_{i,\delta}$.
The AS preconditioning writes: $M_{AS,\delta}^{-1} = \displaystyle\sum_{i=1}^{p}R_{i,\delta}^{T}A_{i,\delta}^{-1}R_{i,\delta} $.

Cai \& Sarkis \cite{RAS} introduced the restriction matrix  $\tilde{R}_{i,\delta}$ on a non-overlapping sub-domain $W_{i,0}$, and then derived the Restricted Additive Schwarz (RAS)  iterative process as:
\begin{equation}
u^{k} = u^{k-1} + M_{RAS,\delta}^{-1}\left(f-Au^{k-1}\right), \, \textrm{with} \, M_{RAS,\delta}^{-1} = \displaystyle\sum_{i=1}^{p}\tilde{R}_{i,\delta}^{T}A_{i,\delta}^{-1}R_{i,\delta} \label{RAS}
\end{equation}
They showed experimentally that the RAS exhibits a faster convergence than the AS, as Efstathiou \& Gander demonstrated in \cite{MR2058877} for the Poisson problem, leading to a better preconditioning that depends on the number of sub-domains. We note, that Cai \& al \cite{MR2116922} develop extensions of RAS for symmetric positive definite problems using the so-called harmonic overlaps (RASHO).

When it is applied to linear problems, the RAS has a pure linear rate of convergence / divergence. When it converges, its convergence can be enhanced with optimized boundary conditions giving the ORAS method of St Cyr \& al \cite{MR2357620}. In this case, the transmission condition in GSAM takes the form $\Lambda_i$ to be the normal derivative ( Neumann boundary condition). Then, an optimisation problem is done to minimize the amplification factor of the Schwarz method with this Robin coefficient in the  Fourier space. The drawback of this method is that it can only be applied to separable operators, and need regular step size and periodic boundary conditions in the direction orthogonal to the interface to be mathematically valid. Nevertheless, if it is not the case, the parameters in the Robin conditions are set based on this postulate and applied in the current case.

This Neumann-Dirichlet map is related to the Schur complement of the discrete operator (see for example Natarajan \cite{MR1453564} or Steinbach \cite{MR2149358}).

Saad \& Li \cite{SchurRAS} introduced the SchurRAS method  based on the ILU factorisation of the local operators $A_{i,\delta}$ present in the RAS method. 
Magoules \& al  \cite{MR2356897} introduced the  patch substructuring methods and demonstrated its  equivalence with the overlapping Schwarz methods. In this work  the Dirichlet 
and Neumann boundary conditions present in the Schwarz alternated method have been replaced by Robin boundary conditions to enhance the convergence rate. 
The patch method consists in introducing an overlap in the Schur complement technique. These two techniques  take care of the data locality in order to avoid global communications involving all sub-domains. Nevertheless, as in the RAS preconditioning technique,  the  main drawback of this locality  is a decreasing of the preconditioning efficiency with respect to the number of sub-domains. 

Another related work that takes care to involve all  the sub-domains present in the Schur complement is the substructuring method with a suitable preconditioner for the reduced equation  of Bramble \& al \cite{MR842125}, Carvalho \& al \cite{MR1856298} \cite{MR1831933}, Khoromskij \& Wittum \cite{MR2045003}. Let us describe the iterative substructuring method and the preconditioning by the additive Schwarz preconditioner for the Schur complement  reduced equation on the interface problem designed  by \cite{MR1856298}. Let $\Gamma$ be the set of all the indices of the mesh points which belong to the interfaces between the sub-domains. Grouping together the unknowns associated to points of the mesh corresponding to $\Gamma$ into the vector $u_\Gamma$ and the ones corresponding to the other unknowns (corresponding to the points of mesh associated to the interior I of  sub-domains)  into the vector $u_I$, we get the reordered problem:
\begin{eqnarray} \begin{pmatrix} A_{II} & A_{I\Gamma} \\ A_{\Gamma I} & A_{\Gamma} \end{pmatrix} \begin{pmatrix} u_I \\ u_\Gamma \end{pmatrix} &=& \begin{pmatrix} f_I \\ f_\Gamma \end{pmatrix} \label{Schur}
\end{eqnarray}
Eliminating the unknowns $u_I$ from the second block row of \eqref{Schur} leads to the following reduced equation for $u_\Gamma$:
\begin{eqnarray}
S u_\Gamma = f_\Gamma - A_{\Gamma I} A_{II}^{-1} f_I,
\label{eq:Schur1}
\end{eqnarray}
where
\begin{eqnarray}
S&=& A_{\Gamma \Gamma} - A_{\Gamma I}A_{II}^{-1}A_{I \Gamma}
\label{eq:Schur2}
\end{eqnarray}
is the Schur complement of the matrix $A_{II}$ in $A$.
Let be  $\Gamma_i = \partial \Omega_i \backslash \partial \Omega$. Let $R_{\Gamma_i} : \Gamma \rightarrow \Gamma_i$  be the canonical pointwise restriction which maps vectors on $\Gamma$ into defined vectors on $\Gamma_i$, and let be $R_{\Gamma_i}^T:\Gamma_i \rightarrow \Gamma$  its transposed. The Schur complement matrix \eqref{eq:Schur2} can also be written as:
\begin{eqnarray}
S=\sum_{i=1}^p R_{\Gamma_i}^T S^{(i)} R_{\Gamma_i}
\label{eq:Schur3}
\end{eqnarray}
where
\begin{eqnarray}
S^{(i)}&=& A_{\Gamma_i \Gamma_i}^{(i)} - A_{\Gamma_i i}A_{ii}^{-1}A_{i \Gamma_i}
\label{eq:Schur4}
\end{eqnarray}
is referred to the local Schur complement associated with the sub-domain $\Omega_i$. $S^{(i)}$ that involves the
submatrices from the local  matrix $A^{(i)}$ which is defined as
\begin{eqnarray} A^{(i)}&=&\begin{pmatrix} A_{ii} & A_{i \Gamma_i} \\ A_{\Gamma_i i} & A_{\Gamma_i \Gamma_i} \end{pmatrix}
\label{eq:Schur5}
\end{eqnarray}

Then they defined a BPS (Bramble, Pasciak \& Schatz \cite{MR842125}) preconditioner which is based on the set $V$ which gathers the cross  points between sub-domains (i.e points that belong to more than two sub-domains) and  the sets $E_i$ of interface points 
(without the cross points in $V$) 
\begin{eqnarray}
E_i = (\partial \Omega_j \cap \partial \Omega_l)-V \\
\Gamma =(\bigcup_{i=1}^m E_i) \cup V
\end{eqnarray}
The operator $R_i$ defined the standard  pointwise restriction of nodal values on $E_i$  while operator $R_V$ defined the canonical restriction on $V$. Then a coarse mesh is associated with the sub-domains  and 
an interpolation operator $R^T$ is defined. This operator corresponds to the linear interpolation between two adjacent cross points $V_j$ $V_i$ in order to define values on the edge $E_i$ . This allows to define  $A_H$ the Galerkin cross grid operator $A_H= R A R^T$.  They deduced a very close variant of BPS preconditioner that can be written as:
\begin{eqnarray}
M_{BPS} = \sum_{E_i} R^T_i S_{ii} R_i + R^T A_H^{-1} R
\label{eq:Giraud1}
\end{eqnarray}
  They defined a coarse-space operator 
\begin{eqnarray}
\Lambda_0 = R_0 S R_0^T
\label{eq:Giraud2}
\end{eqnarray}
where $R_0:U \rightarrow U_0$ is a restriction operator which maps full vector of $U$ into vector in $U_0$ where $U_0$ is a q-dimensional subspace of $U$ the algebraical space of nodal vectors where the Schur complement matrix is defined.
\begin{eqnarray}
M_{BPS} = \sum_{E_i} R^T_i \tilde{S}_{ii} R_i + R_0^T \Lambda_0^{-1} R_0
\label{eq:Giraud3}
\end{eqnarray}
where $\tilde{S}_{ii}$ is an approximation of $S_{ii}$. The definition of $U_0$ gives different preconditioners: Vertex-based coarse space, sub-domain-based coarse space, edge-based coarse space,
depending on the set of points of the interface $\Gamma$  that are involved.
From the implementation practical point of view, the coarse matrix $\Lambda_0$ is constructed once and   involve  matrix vector  products of the local Schur complement only.
\begin{enumerate}
\item The advantages of this method, is to defined the two-level preconditioner only on the interface $\Gamma$. It is intimately related to the Schur complement operator defined on the interface. 
\item The drawback is to have to define a priori the coarse space $U_0$ without any knowledge of the solution behavior. Consequently it can be expensive in term of number of coarse space vectors, specifically for 3D problems where cross-points between sub-domains in 2D, become cross-regions between sub-domains in 3D .
\end{enumerate}

Our approach will follow the same spirit as this two-level preconditioning working only on the interface. But we still work on the system $Ax=b$ and not $Su_\Gamma = g_\Gamma$ and we use an \emph{a posteriori} knowledge of the global Dirichlet to Neumann map that is  based on the pure linear convergence/divergence of the RAS to define the coarse space (equivalent of the definition of $U_0$).

~\\
The plan of this paper is the following. Section \ref{section-ARAS} will derive the Aitken-Schwarz preconditioning. Section \ref{section-approx-explicit} will focus on coarse approximation of the solution at the artificial interfaces, notably with a random set of orthogonal vectors and an orthogonal set of vectors obtained through the SVD of the Schwarz interface solutions. Then, section \ref{section-convergence} proposes a study of convergence of the ARAS class preconditioners. Eventually  numerical tests are provided on academic problems in section \ref{section-results-academic} and industrial problems in section \ref{section-results-ind}.

\section{Aitken-Schwarz method derived as preconditioning technique}
\label{section-ARAS}

In this section we study the integration of the Aitken's acceleration into a Richardson process in order to formulate a preconditioning technique based on Aitken. More precisely, we propose to enhance the RAS preconditioning technique, presented in section \ref{section-introduction}, by the Aitken's acceleration. We first present the mechanism of the method and develop the equation to extract a corresponding Richardson's process. Then we point out that the method in its simple form does not exhibit the complete acceleration after one application and need an update as when the method is used as solver. The result is a multiplicative preconditioner based on the Aitken RAS preconditioner. Finally we present those preconditioners in their approximated form in order to save computing.

		\subsection{The Aitken Restricted Additive Schwarz preconditioner: ARAS}
		\label{subsection-ARAS}
Let $\Gamma_i= (I_{m_{i,\delta}}-{R}_{i,\delta}^{T}) W_{i,\delta}$ be the interface associated to $W_{i,\delta}$ and  $\Gamma = \cup_{i=1}^p \Gamma_{i}$ be the global interface.
Then $u_{|\Gamma} \in \mathbb{R}^{n}$ is the restriction of the solution $u \in \mathbb{R}^{m}$ on the $\Gamma$ interface and $e_{|\Gamma}^k= u_{|\Gamma}^k -u_{|\Gamma}^{\infty}$  is the error of an iteration of a RAS iterative process, equation \eqref{eq:RAS} at the interface $\Gamma$.
\begin{equation}\label{eq:RAS}
u^{k}= u^{k-1}+M_{RAS,\delta}^{-1}(f-Au^{k-1})
\end{equation}

In section \ref{section-introduction} we wrote that the Schwarz iterative method has a pure linear convergence. This property enables us to use the Aitken's technique presented the same section. The previously mentioned $\cal{G}$ iterative process is replaced by a RAS iterative process. 

Using the linear convergence property of the RAS method, we would like to write a preconditioner which includes the Aitken's acceleration process.
We introduce a restriction operator $R_{\Gamma} \in \mathbb{R}^{n\times m}$ from $W$ to the global artificial interface $\Gamma$, with $R_{\Gamma}R_{\Gamma}^{T}=I_{n}$.
The Aitken Restricted Additive Schwarz (ARAS) must generate a sequence of solution on the interface $\Gamma$, and accelerate the convergence of the Schwarz process from this original sequence. Then the accelerated solution on the interface replaces the last one. This could be written combining an AS or RAS process eq.\eqref{subeq-1:ARAS}) with the Aitken process written in $\mathbb{R}^{m\times m}$ eq.\eqref{subeq-2:ARAS} and subtracting the Schwarz solution which is not extrapolated on $\Gamma$ eq.\eqref{subeq-3:ARAS}. We can write the following approximation $u^{*}$ of the solution $u$: 
\begin{subequations}\label{eq:ARAS}
\begin{flalign}
u^{*}= & \mbox{ }u^{k-1}+M_{RAS,\delta}^{-1}(f-Au^{k-1}) \label{subeq-1:ARAS}\\
& + R_{\Gamma}^{T}\left(I_{n}-P\right)^{-1}\left(u_{|\Gamma}^k-Pu_{|\Gamma}^{k-1}\right)\label{subeq-2:ARAS}\\
& - R_{\Gamma}^{T}I_{n}R_{\Gamma}\left(u^{k-1}+M_{RAS,\delta}^{-1}(f-Au^{k-1})\right)\label{subeq-3:ARAS}
\end{flalign}
\end{subequations}

We would like to write $u^{*}$ as an iterated solution derived from an iterative process of the form $u^{*}= u^{k-1} + M_{ARAS,\delta}^{-1}\left(f-Au^{k-1}\right)$, where $M_{ARAS,\delta}^{-1}$ is the Aitken-RAS preconditioner.

First of all, we write an expression of eq.\eqref{subeq-2:ARAS} depending on eq.\eqref{eq:RAS} and which only involves the iterated solution $u^{k-1} \in \mathbb{R}^{m}$, as follows: 
\begin{eqnarray*}
eq.\eqref{subeq-2:ARAS} & := &R_{\Gamma}^{T}\left(I_{n}-P\right)^{-1}\left(u_{|\Gamma}^k-Pu_{|\Gamma}^{k-1}\right)\\
& = & R_{\Gamma}^{T}\left(I_{n}-P\right)^{-1}R_{\Gamma} \left(R_{\Gamma}^{T}I_{n}R_{\Gamma}u^{k}-R_{\Gamma}^{T}PR_{\Gamma}u^{k-1} \right)\\
& \downarrow & \mbox{ with eq.\eqref{eq:RAS} }\\
& = & R_{\Gamma}^{T}\left(I_{n}-P\right)^{-1}R_{\Gamma} \left(R_{\Gamma}^{T}I_{n}R_{\Gamma}\left( u^{k-1} + M_{RAS,\delta}^{-1}\left(f-Au^{k-1}\right)\right)\right.\\& & \left. - R_{\Gamma}^{T}PR_{\Gamma}u^{k-1}\right)\\
& = & R_{\Gamma}^{T}\left(I_{n}-P\right)^{-1}R_{\Gamma}R_{\Gamma}^{T}I_{n}R_{\Gamma}\left( u^{k-1} + M_{RAS,\delta}^{-1}\left(f-Au^{k-1}\right)\right)\\& &- R_{\Gamma}^{T}\left(I_{n}-P\right)^{-1}R_{\Gamma}R_{\Gamma}^{T}PR_{\Gamma}u^{k-1}\\
& = & R_{\Gamma}^{T}\left(I_{n}-P\right)^{-1}R_{\Gamma}\left( u^{k-1} + M_{RAS,\delta}^{-1}\left(f-Au^{k-1}\right)\right)\\& &- R_{\Gamma}^{T}\left(I_{n}-P\right)^{-1}PR_{\Gamma}u^{k-1}
\end{eqnarray*}

Then, we re-write eq.\eqref{eq:ARAS} with this new expression of eq.\eqref{subeq-2:ARAS} as follows:
\begin{eqnarray*}
u^{*}& = & u^{k-1}+M_{RAS,\delta}^{-1}(f-Au^{k-1})\\& & + R_{\Gamma}^{T}\left(I_{n}-P\right)^{-1}R_{\Gamma}\left( u^{k-1} + M_{RAS,\delta}^{-1}\left(f-Au^{k-1}\right)\right)\\ &&- R_{\Gamma}^{T}\left(I_{n}-P\right)^{-1}PR_{\Gamma}u^{k-1} - R_{\Gamma}^{T}I_{n}R_{\Gamma}\left(u^{k-1}+M_{RAS,\delta}^{-1}(f-Au^{k-1})\right)\\
& \downarrow & \mbox{ factorizing by }\left( u^{k-1} + M_{RAS,\delta}^{-1}\left(f-Au^{k-1}\right)\right)\\
& = & \left(I_m - R_{\Gamma}^{T}I_{n}R_{\Gamma} + R_{\Gamma}^{T}\left(I_{n}-P\right)^{-1}R_{\Gamma} \right)\left(u^{k-1} + M_{RAS,\delta}^{-1}\left(f-Au^{k-1}\right)\right)\\
& &- R_{\Gamma}^{T}\left(I_{n}-P\right)^{-1}PR_{\Gamma}u^{k-1}\\
& \downarrow & \mbox{ isolating } u^{k-1} \mbox{ from } M_{RAS,\delta}^{-1}\left(f-Au^{k-1}\right)\\
& = & u^{k-1} + \left(-R_{\Gamma}^{T}I_{n}R_{\Gamma}+R_{\Gamma}^{T}\left(I_{n}-P\right)^{-1}R_{\Gamma}-R_{\Gamma}^{T}\left(I_{n}-P\right)^{-1}PR_{\Gamma}\right)u^{k-1} \\
& & + \left(I_m - R_{\Gamma}^{T}I_{n}R_{\Gamma} + R_{\Gamma}^{T}\left(I_{n}-P\right)^{-1}R_{\Gamma} \right)M_{RAS,\delta}^{-1}\left(f-Au^{k-1}\right)
\end{eqnarray*}

One can simplify $E = \left( R_{\Gamma}^{T}\left(I_{n}-P\right)^{-1}R_{\Gamma}-R_{\Gamma}^{T}\left(I_{n}-P\right)^{-1}PR_{\Gamma}\right)$ as follows:
\begin{eqnarray*}
E & = & R_{\Gamma}^{T}\left(I_{n}-P\right)^{-1}R_{\Gamma}\left(R_{\Gamma}^{T}I_{n}R_{\Gamma} - R_{\Gamma}^{T}PR_{\Gamma}\right)\\
& = & R_{\Gamma}^{T}\left(I_{n}-P\right)^{-1}R_{\Gamma}R_{\Gamma}^{T}\left(I_{n}-P\right)R_{\Gamma}\\
& = & R_{\Gamma}^{T}\left(I_{n}-P\right)^{-1}\left(I_{n}-P\right)R_{\Gamma}\\
& = & R_{\Gamma}^{T}I_{n}R_{\Gamma}\\
\end{eqnarray*}

And then writes,
\begin{eqnarray*}
u^{*}&=&u^{k-1} + \left(-R_{\Gamma}^{T}I_{n}R_{\Gamma}+R_{\Gamma}^{T}I_{n}R_{\Gamma}\right)u^{k-1} \\
& & + \left(I_m - R_{\Gamma}^{T}I_{n}R_{\Gamma} + R_{\Gamma}^{T}\left(I_{n}-P\right)^{-1}R_{\Gamma} \right)M_{RAS,\delta}^{-1}\left(f-Au^{k-1}\right)\\
& = &  u^{k-1} + \left(I_m - R_{\Gamma}^{T}I_{n}R_{\Gamma} + R_{\Gamma}^{T}\left(I_{n}-P\right)^{-1}R_{\Gamma} \right)M_{RAS,\delta}^{-1}\left(f-Au^{k-1}\right)\\
& = &  u^{k-1} + \left(I_m + R_{\Gamma}^{T}\left(\left(I_{n}-P\right)^{-1}-I_{n}\right)R_{\Gamma}\right)M_{RAS,\delta}^{-1}\left(f-Au^{k-1}\right)
\end{eqnarray*}

Hence the formulation eq.\eqref{eq:ARAS} leads to an expression of an iterated solution $u^{*}$:
\begin{equation*}
u^{*}= u^{k-1} + \left(I_{m}+R_{\Gamma}^{T}\left( \left(I_{n}-P\right)^{-1}-I_{n}\right)R_{\Gamma}\right)M_{RAS,\delta}^{-1}\left(f-Au^{k-1}\right)
\end{equation*}

This iterated solution $u^{*}$ can be seen as an accelerated solution of the RAS iterative process. Drawing our inspiration from the Stephensen's method \cite{bookStoerBulirsch}, we build a new sequence of iterates from the solutions accelerated by the Aitken's acceleration method. Then, one considers $u^{*}$ as a new $u^{k}$ and writes the following ARAS iterative process:
\begin{equation}\label{ARAS-process}
u^{k}= u^{k-1} + \left(I_{m}+R_{\Gamma}^{T}\left( \left(I_{n}-P\right)^{-1}-I_{n}\right)R_{\Gamma}\right)M_{RAS,\delta}^{-1}\left(f-Au^{k-1}\right)
\end{equation}

Then we defined the ARAS preconditioner as
\begin{equation}\label{ARAS-preconditioner}
M_{ARAS,\delta}^{-1} = \left(I_{m}+R_{\Gamma}^{T}\left( \left(I_{n}-P\right)^{-1}-I_{n}\right)R_{\Gamma}\right)\displaystyle\sum_{i=1}^{p}\tilde{R}_{i,\delta}^{T}A_{i,\delta}^{-1}R_{i,\delta}
\end{equation}

\begin{remark}
The ARAS preconditioner can be considered as a two-level additive preconditioner. The preconditioner consists in computing a solution on an entire domain applying the RAS preconditioner and add components computed only on the interface $\Gamma$. 
\end{remark}

		\subsection{Composite Multiplicative form of ARAS: ARAS2} 
		\label{subsection-ARAS2}
If $P$ is known exactly, the ARAS process written in the equation \eqref{ARAS-process} needs two steps to converge to the 
solution $u$ with an initial guess $u^0=0$. Then we have:

\begin{proposition}\label{prop-Am1}
If $P$ is known exactly then we have $$A^{-1}=\left(2M_{ARAS,\delta}^{-1}-M_{ARAS,\delta}^{-1}AM_{ARAS,\delta}^{-1}\right)$$  that leads  $\left(I-M_{ARAS,\delta}^{-1}A\right)$ to be a nilpotent matrix of degree 2.
~\\ \noindent{\bf Proof} We consider the postulate: "If $P$ is known exactly, the ARAS process written in eq.\eqref{ARAS-process} needs two steps to converge to the solution with an initial guess $u^0=0$".

We write the two first iterations of the ARAS process for any initial guess $u^0 \in \mathbb{R}^m$ and for all $f\in\mathbb{R}^m$:
\begin{eqnarray*}
u^{1} & = & u^{0}+M_{ARAS,\delta}^{-1}\left(f-Au^{0}\right)
\end{eqnarray*}

And the second iterations leads to:
\begin{eqnarray*}
u^{2} & = & u^{1}+M_{ARAS,\delta}^{-1}\left(f-Au^{1}\right)\\
& = & u^{0}+M_{ARAS,\delta}^{-1}\left(f-Au^{0}\right) + M_{ARAS,\delta}^{-1}\left(f-A \left(u^{0}+M_{ARAS,\delta}^{-1}\left(f-Au^{0}\right) \right) \right)
\end{eqnarray*}

Let $u^0=0$, then,
\begin{eqnarray*}
u^{2} & = & M_{ARAS,\delta}^{-1}f + M_{ARAS,\delta}^{-1}\left(f-A \left(M_{ARAS,\delta}^{-1}f \right) \right)\\
& = & \left(2M_{ARAS,\delta}^{-1} - M_{ARAS,\delta}^{-1}AM_{ARAS,\delta}^{-1}\right)f\\
& = & u
\end{eqnarray*}

Since this expression is true for all $f\in\mathbb{R}^m$ we can write:
\begin{eqnarray*}
A^{-1} & = & 2M_{ARAS,\delta}^{-1} - M_{ARAS,\delta}^{-1}AM_{ARAS,\delta}^{-1}\\
\end{eqnarray*}

Now we can write: 
\begin{eqnarray*}
u & = & \left(2M_{ARAS,\delta}^{-1} - M_{ARAS,\delta}^{-1}AM_{ARAS,\delta}^{-1}\right)f\\
  & = & \left(I_m + I_m - M_{ARAS,\delta}^{-1}A\right)M_{ARAS,\delta}^{-1}f\\
  & = & M_{ARAS,\delta}^{-1}f + \left(I_m - M_{ARAS,\delta}^{-1}A\right)M_{ARAS,\delta}^{-1}f\\
  & \downarrow & \mbox{ with } Au = f\\
  & = & M_{ARAS,\delta}^{-1}Au + \left(I_m - M_{ARAS,\delta}^{-1}A\right)M_{ARAS,\delta}^{-1}Au\\
\end{eqnarray*}

Thus,
\begin{eqnarray*}
\left(I_m - M_{ARAS,\delta}^{-1}A\right)u & = & \left(I_m - M_{ARAS,\delta}^{-1}A\right)M_{ARAS,\delta}^{-1}Au\\
\end{eqnarray*}

Which is equivalent to 
\begin{eqnarray*}
0 & = &\left(I_m - M_{ARAS,\delta}^{-1}A\right)^{2}u \mbox{, } \forall u\in \mathbb{R}^{m} \\
\end{eqnarray*}

Hence $\left(I_m - M_{ARAS,\delta}^{-1}A\right)$ is nilpotent of degree $2$. $ [] $
\end{proposition}

The previous proposition leads to an approximation of $A^{-1}$ written from the $2$ first iterations of the ARAS iterative process \eqref{ARAS-process}.
Those $2$ iterations compute the Schwarz solutions sequence on the interface needed in order to accelerate the Schwarz method by the Aitken's acceleration.  We now write $2$ iterations of the ARAS iterative process \eqref{ARAS-process} for any initial guess and for all $u^{k-1} \in \mathbb{R}^{m}$.

\begin{eqnarray*}
u^{k+1} & = & u^{k-1} +M_{ARAS,\delta}^{-1}\left(f-Au^{k-1}\right) \\&&+M_{ARAS,\delta}^{-1}\left(f-A\left(u^{k-1}+M_{ARAS,\delta}^{-1}\left(f-Au^{k-1}\right)\right)\right) \\
& = &  u^{k-1} +M_{ARAS,\delta}^{-1}f-M_{ARAS,\delta}^{-1}Au^{k-1} \\&&+M_{ARAS,\delta}^{-1}f-M_{ARAS,\delta}^{-1}Au^{k-1} - M_{ARAS,\delta}^{-1}AM_{ARAS,\delta}^{-1}\left(f-Au^{k-1}\right) \\
& = &  u^{k-1} +2M_{ARAS,\delta}^{-1}\left(f-Au^{k-1}\right) - M_{ARAS,\delta}^{-1}AM_{ARAS,\delta}^{-1}\left(f-Au^{k-1}\right) \\
& = &  u^{k-1} +\left(2M_{ARAS,\delta}^{-1} - M_{ARAS,\delta}^{-1}AM_{ARAS,\delta}^{-1}\right)\left(f-Au^{k-1}\right) \\
\end{eqnarray*}

Then we defined the ARAS2 preconditioner as
\begin{equation}
M_{ARAS2,\delta}^{-1} = 2M_{ARAS,\delta}^{-1} - M_{ARAS,\delta}^{-1}AM_{ARAS,\delta}^{-1}
\end{equation}

\begin{remark}
According to the linear algebra literature about preconditioning technique \cite{MR2101970}, the ARAS2 preconditioner can be considered as a composite multilevel preconditioner. Actually, the ARAS2 preconditioner is a multiplicative form of ARAS which is itself an additive preconditioner adding an operation on the entire domain with RAS and an operation on a coarse interface with the Aitken formula.
\end{remark}

		\subsection{Approximated form of ARAS and ARAS2}
		\label{subsection-approx-ARAS-ARAS2}
As the previous subsection suggests, since $P$ is known exactly there is no need to use ARAS as a preconditioning technique. Nevertheless, when $P$ is approximated, the Aitken's acceleration of the convergence depends on the local domain solving accuracy, and the cost of the building of an exact $P$ depends on the size $n$. 
This is why $P$ is numerically approximated by $P_{\mathbb{U}_q}$ as in \cite{pareng09SVD}, defining $q \le n$ orthogonal vectors $U_i \in \mathbb{R}^n$, defining the columns of the matrix $\mathbb{U}_q \in R^{n\times q}$. Then it makes sense to use $P_{\mathbb{U}_q}$ in the ARAS preconditioning technique to define the ARAS(q) preconditioner:
\begin{equation}
M^{-1}_{ARAS(q),\delta} = \left(I_{m}+R_{\Gamma}^{T}\mathbb{U}_q\left( \left(I_{q}-{P_{\mathbb{U}_q}}\right)^{-1}-I_{q}\right)\mathbb{U}_q^T R_{\Gamma}\right)\displaystyle\sum_{i=1}^{p}\tilde{R}_{i,\delta}^{T}A_{i,\delta}^{-1}R_{i,\delta} 
\end{equation}
and ARAS2(q),
\begin{equation}
M^{-1}_{ARAS2(q),\delta} = 2M_{ARAS(q),\delta}^{-1} - M_{ARAS(q),\delta}^{-1}AM_{ARAS(q),\delta}^{-1}
\end{equation}

Different kind of approximation techniques of the error transfer operator matrix $P$ have been proposed in the work of  \cite{Garbey1,mg_dtd_fd,barberou,Baranger,Frullone}. Nevertheless, we are interested in algebraic ways to compute an Aitken acceleration. In section \ref{section-approx-explicit} we proposed two fully algebraic approaches.

\section{Basis to approximate the interface solution}
\label{section-approx-explicit}

\begin{figure}[h]
\begin{center}
\includegraphics[scale=0.8]{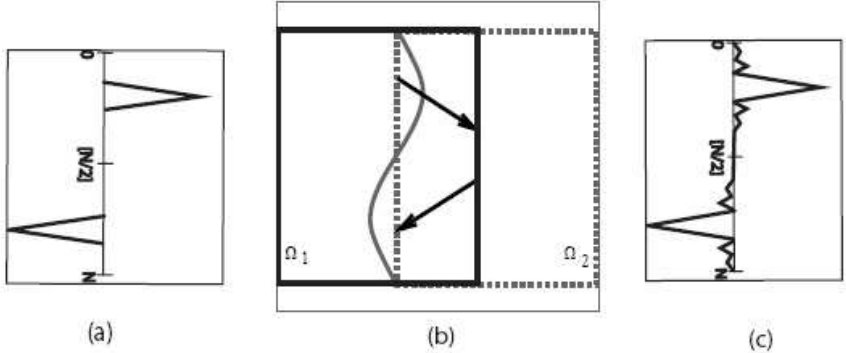}
\end{center}
 \caption{Steps to build the $\hat{P}$ matrix}
\label{figure-explicit-construct}
\end{figure}

We focus here on an algebraic way to compute an Aitken acceleration of a sequence of Schwarz solutions on the interface. The global approach consists on an explicit building of $\hat{P}$ computing how the spanning vectors $U_i$ are modified by the Schwarz iteration. Figure \ref{figure-explicit-construct} illustrates the steps for constructing the matrix $\hat{P}$. Step (a) starts from the spanning vector on the interface $R_{\Gamma}^{T}U_i$ and gets its value on the interface in the physical space. Then step (b) performs a complete Schwarz iteration with zero local right hand sides and homogeneous boundary conditions on the others artificial interfaces. Step (c) decomposes the trace solution on the interface in the spanning vector set $\mathbb{U}_q$. Thus, we obtain the column $i$ of the matrix $\hat{P}$.

The full computation of $\hat{P}$ can be done in parallel, but it needs as much local domain solution as the number of interface points (i.e the size of the matrix $\hat{P}$). 

Its adaptive computation is required to save computing. This methodology was first used with Fourier basis functions \cite{Frullone,dufaud-parcfd09}. 
This section focuses on the definition of orthogonal "base" $\mathbb{U}_q$ that will extend this adaptive computation in a general context. In the following, we use the term "base" to denote a spanning vectors set that defines the approximation space. The key point of these preconditioners's efficiency  is the choice of this orthogonal "base" $\mathbb{U}_q$. It must be sufficiently rich to numerically represent the solution at the interface, but it has to be not too large for the computation's efficiency. 

We first propose a "naive" approach consisting of selecting an arbitrary set of orthogonalized random vectors to generate the space to approximate the solution. Secondly, we represent the solution in a space arising from the singular value decomposition of a sequence of Schwarz solutions. Doing this, we assume to represent the main modes of the solutions.

		\subsection{Orthogonal "base" arising from an arbitrary coarse algebraic approximation of the interface.}
		\label{subsection-random-base}
A choice consists in having a coarse representation of the interface's solution $u \in \mathbb{R}^n$ from an algebraical point of view. 
Nevertheless, it is not possible to take a subset of $q$ vectors of the canonical base of $\mathbb{R}^n$, as if some components of $u$ are not reachable by the "base" $\mathbb{U}_q$, then the approximation $||u-\mathbb{U}_q ( \mathbb{U}_q^t u)||$ will be very bad. This reason leads us to define $\mathbb{U}_q$ as a set of orthogonal vectors where each component is coming from a random process in order that each vector can contribute to a part of the searched solution at the interface. 
We split $q$ such as $q=\displaystyle\sum_{i=1}^{p}q_{i}$ and we associate $q_i$ random vectors to the interface $\Gamma_i$, $1 \le i \le p$. Then these $q_i$ vectors are orthogonalized to form $q_i$ columns of the "base" $\mathbb{U}_q$.  

This strategy is hazardous but can be a simple way to improve the convergence of a Schwarz process without knowledge of the problem and the mesh. 

The orthogonal "base" $\mathbb{U}_q$ is obtained applying the same principle as illustrated in Figure \ref{figure-explicit-construct}, leading to Algorithm \ref{algorithm-random}.
\begin{algorithm}[ht]
\caption{Vectorial Aitken's acceleration in an arbitrary built space without inversion}
\label{algorithm-random}
\begin{algorithmic}[1]
\REQUIRE {${\cal{G}}:   \mathbb{R}^n \rightarrow \mathbb{R}^n$ an iterative method having a pure linear convergence}
\STATE {Compute $q$ random vectors $v_i \in \mathbb{R}^n$  following a uniform law on $ \left[ 0, 1\right] $}
\STATE {Orthogonalize those $q$ vectors to form $\mathbb{U}_q \in \mathbb{R}^{n \times q}$}
\STATE {Apply {\bf one} iterate of $\cal{G}$ on homogeneous problem,  $\mathbb{U}_q  \rightarrow W ={\cal{G}}(\mathbb{U}_q)$ }
\STATE {Set  $\hat{P}=\mathbb{U}_q^t W$}
\STATE { $\hat{\xi} = (I_q - \hat{P})^{-1}\;(\hat{u}^1-\hat{P}\;\hat{u}^0)$} \COMMENT{Aitken Formula}
\STATE {$\xi=\mathbb{U}_q \hat{\xi}$}
\end{algorithmic}
\end{algorithm}

The lack of this method is that there is no possibility to control the quality of the base to perform the acceleration. A more controllable method will be preferred. In the following subsection, we propose a different starting point to build the base. The main idea will be in the fact that we can compress the vectorial sequence using  a Singular Value Decomposition. Since the $\mathbb{U}_q$ matrix is built, $P$ is built the same way.

	\subsection{Approximation compressing the vectorial sequence}
	\label{subsection-SVD-base}

A totally algebraic method based on the Singular Value Decomposition of the Schwarz solutions on the interface has been proposed when the modes of the error could be strongly coupled \cite{pareng09SVD}. This method offers the possibility for the Aitken Schwarz method to be used on a large class of problem without mesh consideration.
Moreover, when computing an Aitken acceleration, the main difficulty is to invert the matrix $\left[ E^{n-1},\ldots,E^0 \right]$ which can be close to singular.
In a computation, most of the time is consumed solving some noise that does not actually contribute to the solution. The singular value decomposition offer a tools to concentrate the effort only on the main parts of the solution.

		\subsubsection{The singular value decomposition}
		\label{subsubsection-SVD}
A singular-value decomposition (SVD) of a real $n \times m$ $(n > m)$ matrix $A$ is its factorization into the product of three matrices:
\begin{eqnarray}
A&=&\mathbb{U} \Sigma \mathbb{V}^*,
\end{eqnarray}
where $\mathbb{U}=\left[U_1,\ldots, U_m\right]$ is a $n \times m$ matrix with orthonormal columns, $\Sigma$ is a $n \times m$  non-negative diagonal matrix  with $\Sigma_{ii}=\sigma_i, \, 1\leq i \leq m$ and the $m \times m$ matrix $\mathbb{V}=\left[ V_1, \ldots, V_m \right]$ is orthogonal. The left $\mathbb{U}$ and right $\mathbb{V}$ singular vectors are the eigenvectors of $AA^*$ and $A^*A$ respectively.  
It readily follows that $Av_i=\sigma_i u_i, \, 1 \leq i \leq m$

We are going to recall some properties of the SVD. Assume that the $\sigma_i, 1\leq i \leq m$ are ordered in decreasing order and there exists $r$ such that $\sigma_r>0$ while $\sigma_r+1=0$. Then $A$ can be decomposed  in a dyadic decomposition: 
\begin{eqnarray} 
A=\sigma_1 U_1 V_1^*+\sigma_2 U_2 V_2^*+ \ldots + \sigma_r U_r V_r^*. \label{dyadic}
\end{eqnarray}
This means that SVD provides a way to find optimal lower dimensional approximations of a given series of data. More precisely, it produces an orthonormal base for representing the data series in a certain least squares optimal sense. This can be summarized  by the theorem of Schmidt-Eckart-Young-Mirsky: 
\begin{theorem} \label{thm1}
A non unique minimizer $X_*$ of the problem $\min_{X, rank X=k} || A - X||_2= \sigma_{k+1}(A)$, provided that $\sigma_k > \sigma_{k+1}$, is obtained by truncating the dyadic decomposition  of \eqref{dyadic} to contain the first $k$ terms:
$X_*=\sigma_1 U_1 V_1^*+\sigma_2 U_2 V_2^*+ \ldots + \sigma_k U_k V_k^*$
\end{theorem}
The SVD of a matrix is well conditioned with respect to perturbations of its entries. Consider the matrix $A,B \in \mathbb{R}^n$,  the Fan inequalities write $\sigma_{r+s+1}(A+B) \leq \sigma_{r+1}(A)+\sigma_{s+1}(B)$ with $r,s \geq 0, \, r+s+1 \leq n$. Considering the perturbation matrix $E$ such that $||E||=O(\epsilon)$, then $|\sigma_i(A+E)-\sigma_i(A)| \leq \sigma_1(E)=||E||_2, \, \forall i$. This property does not hold for eigenvalues decomposition where small perturbations in the matrix entries can cause a large change in the eigenvalues.

This property allows us to search the acceleration of the convergence of the sequence of vectors in the base linked to its SVD.  

\begin{proposition}\label{prop:define-P-SVD}
Let $(u^k)_{1 \leq k \leq q}$ $q$ successive iterates satisfying the pure linear convergence property: $u^k-u^\infty=P (u^{k-1}-u^\infty)$.
Then there exists an orthogonal base $\mathbb{U}_q=\left[ U^1, U^2,\ldots,U^{q} \right]$ of a subset of ${\mathbb R}^n$ such that
\begin{itemize}
\item $\alpha_l^k = \sigma_l V_{kl}^*$ with $(\sigma_l)_{l \in \mathbb{N}}$ decreasing and $|V_{kl}^*| \leq 1 \Rightarrow \forall l \in \left\{1,...,q \right\}$,  \\$|\alpha_l^k| \leq |\sigma_l|$
\item $u^k=\sum_{l=1}^{q} \alpha_l^k U^l, \forall k \in \left\{1,...,q \right\}$
\end{itemize}
One can write:
\begin{eqnarray}
(\alpha_1^{k+1}-\alpha_1^k,\ldots, \alpha_q^{k+1}-\alpha_q^k)^T= \hat{P} (\alpha_1^{k}-\alpha_1^{k-1},\ldots, \alpha_q^{k}-\alpha_q^{k-1})^T
\end{eqnarray}
where $\hat{P}  \stackrel{def}{=}\mathbb{U}_q^* P \mathbb{U}_q$. Moreover $(\alpha^\infty_1,\ldots,\alpha_q^\infty)^T$ obtained by the acceleration process represents the projection of the limit of the sequence \
of vectors in the space generated by $\mathbb{U}_q$.
~\\ \\ \noindent{\bf Proof} By theorem 3.1.3 there exists a SVD decomposition of $\left[u^1,\ldots,u^q\right]=\mathbb{U}_q \Sigma \mathbb{V}^{*}$ and we can identify $\alpha_l^k$ as $\sigma_l V_{kl}^*$. The orth\
onormal property  of $\mathbb{V}$ associated to the decrease of $\sigma_l$ with increasing $l$ leads to have $\alpha_l^k$ bounded by $|\sigma_l|$: $\forall l \in \left\{1,...,q \right\}$,  $|\alpha_l^k| \leq |\sigma_l|$.\\
Taking the pure linear convergence of $u^k$ in the matrix form, and applying $\mathbb{U}_q$ leads to:
\begin{eqnarray}
\mathbb{U}_q^* (u^k-u^\infty)=\mathbb{U}_q^* P \mathbb{U}_q \mathbb{U}_q^*(u^{k-1}-u^\infty) \\
(\alpha_1^k-\gamma_1^\infty,\ldots, \alpha_q^{k}-\gamma_q^\infty)^T= \hat{P} (\alpha_1^{k-1}-\gamma_1^\infty,\ldots, \alpha_q^{k-1}-\gamma_q^\infty)^T
\end{eqnarray}
where $(\gamma_j^\infty)_{1\leq j \leq q}$ represents the projection of $u^\infty$ on the $span\left\{U_1,\ldots,U_q \right\}$. $[]$
\end{proposition}

We can then derive Algorithm \ref{algorithm-SVD-with-inverting}. This algorithm is similar to Algorithm \ref{algorithm-physical-space} since the error transfer operator is defined using the errors of the linear iterative process in a space arising from the Singular Values Decomposition of $q+2$ successive iterates. Therefore the third step of Algorithm \ref{algorithm-physical-space} is equivalent to the sixth step of Algorithm \ref{algorithm-SVD-with-inverting}.
\begin{algorithm}[h]
\caption{Vectorial Aitken's acceleration in the SVD space with inversion}
\label{algorithm-SVD-with-inverting}
\begin{algorithmic}[1]
\REQUIRE {${\cal{G}}:   \mathbb{R}^n \rightarrow \mathbb{R}^n$ an iterative method having a pure linear convergence}
\REQUIRE {$(u^k)_{1 \leq k \leq q+2}$, \; $q+2$ successive iterates of ${\cal{G}}$ to solve the linear system $Au=f$, starting from an arbitrary initial guess $u^0$}
\STATE {Form the SVD decomposition of $Y=\left[ u^{q+2},\ldots,u^1 \right]=\mathbb{U}S V^{T}$}
\STATE {Set $l$ the index such that $l=\max_{1\leq i \leq m+1} \left\{ S(i,i)> tol \right\}$,} \COMMENT{ex.:$tol=10^{-12}$.}
\STATE {Set $\hat{Y}_{1:l,1:l+2}=S_{1:l,1:l}V_{1:l,q-l:q+2}^t$}
\STATE {Set $\hat{E}_{1:l,1:l+1}=\hat{Y}_{1:l,2:l+2}-\hat{Y}_{1:l,1:l+1}$}
\IF{$\hat{E}_{1:l,1:l}$ is non singular}  \STATE{$\hat{P}=\hat{E}_{1:l,2:l+1} \hat{E}_{1:l,1:l}^{-1}$}
\STATE { $\hat{y}^\infty_{1:l,1} = (I_l-\hat{P})^{-1}\;(\hat{Y}_{1:l,l+1}-\hat{P}\hat{Y}_{1:l,l})$} \COMMENT{Aitken Formula}
\STATE {$u^\infty=\mathbb{U}_{:,1:l}\;\hat{y}^\infty_{1:l,1}$}
\ENDIF
\end{algorithmic}
\end{algorithm}

\begin{proposition}\label{prop:algorithm-svd-with-inverting-converges}
Successive applications of Algorithm \ref{algorithm-SVD-with-inverting} converge to the limit $u^\infty$.
~\\ \noindent{\bf Proof} As the sequence of vector $u^k$ converges to a limit $u^\infty$ then we can write $$\Xi=\left[ u^1, \ldots, u^q \right]= \left[ u^\infty, \ldots, u^\infty \right] + E$$ where $E$ is a $n\times q$ matrix with decreasing coefficients with respect to the columns. The SVD of $\Xi^\infty=\left[ u^\infty, \ldots, u^\infty \right]$ leads to have $U^1=u^\infty$ and $\sigma_i(\Xi^\infty)=0,\, i \geq 2$. The fan inequalities lead to have $\sigma_i(\Xi) \leq \sigma_1(E)=||E||_2, i \geq 2$. Consequently, successive applications of Algorithm \ref{algorithm-SVD-with-inverting} decrease the number of non zero singular values at each application. $ [] $
\end{proposition}

In Algorithm \ref{algorithm-SVD-with-inverting} the building of $P$ needs the inversion of the matrix $\hat{E}_{1:l,1:l}$ which can contain very small singular values even if we selected those greater than a certain tolerance. This singular value can deteriorate the ability of $P$ to accelerate the convergence. If it is the case,  we can proceed  inverting this matrix with its SVD, replacing by zeros the singular values less than a tolerance instead of inverting them (see numerical recipes \cite{numerical-recipes}). A more robust algorithm can be obtained without inverting $\hat{E}_{1:l,1:l}$. It consists in building $P$ by applying the iterative method  $\cal{G}$ to the selected columns of $\mathbb{U}_q$ that appears in Algorithm \ref{algorithm-SVD-with-inverting}. Then $\hat{P}=\mathbb{U}^*_{1:n,1:l}{\cal{G}}(\mathbb{U}_{1:n,1:l})$ as done in Algorithm \ref{algorithm-SVD-without-inverting}.

\begin{algorithm}[ht]
\caption{Vectorial Aitken acceleration in the SVD space without inversion}
\label{algorithm-SVD-without-inverting}
\begin{algorithmic}[1]
\REQUIRE {${\cal{G}}:   \mathbb{R}^n \rightarrow \mathbb{R}^n$ an iterative method having a pure linear convergence}
\REQUIRE {$(u^k)_{1 \leq k \leq q+2}$, \; $q+2$ successive iterates of ${\cal{G}}$ to solve the linear system $Au=f$, starting from an arbitrary initial guess $u^0$}
\STATE {Form the SVD decomposition of $Y=\left[ u^{q+2},\ldots,u^1 \right]=\mathbb{U} S V^{T}$}
\STATE {Set the index $l$ such that $l=\max_{1\leq i \leq q+1} \left\{ S(i,i)> tol \right\}$,} \COMMENT{ex.:$tol=10^{-12}$.}
\STATE {Apply {\bf one} iterate of $\cal{G}$ on homogeneous problem, with $l+2$ initial guesses  $\mathbb{U}_{:,1:l}  \rightarrow W_{:,1:l}={\cal{G}}(\mathbb{U}_{:,1:l})$ }
\STATE {Set  $\hat{P}=\mathbb{U}_{:,1:l}^t W_{:,1:l}$}
\STATE {Set $ \hat{Y}_{1:l,1:2}=S_{1:l,1:l} V_{1:l,q+1:q+2}^t$}
\STATE { $\hat{y}^\infty_{1:l,1} = (I_l - \hat{P})^{-1}\;(\hat{Y}_{1:l,2} - \hat{P}\;\hat{Y}_{1:l,1})$} \COMMENT{Aitken Formula}
\STATE {$u^\infty=\mathbb{U}_{:,1:l}\;\hat{y}^\infty_{1:l,1}$}
\end{algorithmic}
\end{algorithm}

\section{Convergence of ARAS and ARAS2 and their approximated form}
	\label{section-convergence}

As an enhancement of the RAS preconditioning technique, ARAS and ARAS2 should have a better convergence rate than the RAS technique. We formulate the convergence rate of a RAS technique considering the linear convergence of the Restricted Additive Schwarz method and extend this formulation to the Aitken's technique. Then we propose a relation between the spectral radius of those methods. 

In the following we note 
\begin{equation}
T_{*} = (I - M_{*}^{-1}A) 
\end{equation}
Any Richardson's process can be written as
\begin{equation}
u^{k} = T_{*} u^{k-1} + c \text{, where } c \in  \mathbb{R}^n \text{is constant}
\end{equation}

\begin{remark}
As the ARAS2 iterative process correspond to $2$ iterations of the ARAS process, we notice that $ T_{ARAS2} = T_{ARAS}^2 $.
\end{remark}

		\subsection{Ideal case}
		\label{subsection-cv-ideal-case}
When building the ARAS preconditioner, we exhibit the fact that $T_{ARAS}$ is nilpotent when the error's transfer operator on the interface $\Gamma$, if $P$ is considered exact. This property gives the following proposition: 
\begin{proposition}\label{prop:rho-if-P-exact}
If $P$ is known exactly then,
\begin{equation}
	\rho\left(T_{ARAS2}\right) = \rho\left(T_{ARAS}\right) = 0
\end{equation}
\noindent{\bf Proof} If $P$ is known exactly then Proposition \ref{prop-Am1} is verified and $T_{ARAS}$ and $T_{ARAS2}$ are nilpotent.
The spectral radius of a nilpotent matrix is equal to $0$. $ [] $
\end{proposition}

\begin{remark}
Obviously, $\rho\left(T_{ARAS2}\right) = \rho\left(T_{ARAS}\right) < \rho\left(T_{RAS}\right)$
\end{remark}

But the matrix $P$ is often numerically computed and then $\rho\left(T_{ARAS}\right)$ is no longer equal to $0$. The value of  $\rho\left(T_{ARAS}\right)$ depends on the accuracy of the local domain solutions and when $P$ is written in another space, depends on the quality of this space.
In the following we propose a framework to study the convergence of $T_{ARAS(q)}$ and $T_{ARAS2(q)}$. The goal is to provide key elements to understand the influence of approximating $P$ in an orthogonal base on the preconditioner. 

		\subsection{Convergence of RAS for an elliptic operator}
		\label{subsection-cv-RAS}

In this subsection we express the convergence rate of a RAS iterative process considering its convergence on the artificial interfaces in proposition \ref{prop:convergence-RAS}. Since we can link the convergence of RAS on the entire domain to the convergence on the interface, it becomes possible to study the effect of modifying the error's transfer operator $P$.

\begin{proposition}\label{prop:convergence-RAS}
Let $A$ be a discretized operator of an elliptic problem on a domain $\Omega$. Let us consider a RAS iterative process such as $T_{RAS} = I - M_{RAS}^{-1}A$ defined on $p$ domains. The data dependencies between domains is located on an artificial interface $\Gamma$.
Then there exists an error's transfer operator on the interface $\Gamma$, $P$ such as there exists a norm $||.||$ for which $||P|| < 1$.
The convergence rate of $T_{RAS}$ is 
\begin{equation}
\rho(T_{RAS}) = \max \left\lbrace |\lambda| \text{ : } \lambda \in \lambda(P) \right\rbrace
\end{equation}
\noindent{\bf Proof} In the case of elliptic problem the maximum principle is observed. Then, for the Schwarz method the error is maximal on artificial interfaces.
We write the error of a Schwarz process starting for the definition of the RAS iterative method:
\begin{equation}
u^{k+1} = T_{RAS}u^{k} + M_{RAS}^{-1}b
\end{equation}	

The convergence of such a process is given by \cite{bookAxelsson} \cite{bookCiarlet}:
\begin{equation}\label{eq:error-process-TRAS}
e^{k} = T_{RAS}^{k}e^{0}
\end{equation}	
On the interface, one can write:
\begin{equation}\label{eq:error-process-P}
e_{|\Gamma}^{k} = P^{k}e_{|\Gamma}^{0}
\end{equation}
The error is maximal on the interface thus, 
\begin{equation}\label{eq:error-RAS-infty}
||e^{k}||_{\infty} = ||e_{|\Gamma}^{k}||_{\infty} 
\end{equation}

Equations (\ref{eq:error-process-TRAS}), (\ref{eq:error-process-P}) and (\ref{eq:error-RAS-infty}) lead to
\begin{equation}
||T_{RAS}^{k}e^{0}||_{\infty} = ||P^{k}e_{|\Gamma}^{0}||_{\infty} < ||P^{k}||_{\infty}||e_{|\Gamma}^{0}||_{\infty} 
\end{equation}
Then we can write,
\begin{eqnarray}
\sup_{||e^{0}||_{\infty}=1}\left( ||T_{RAS}^{k}e^{0}||_{\infty}\right) &=& \sup_{||e_{|\Gamma}^{0}||_{\infty}=1}\left( ||P^{k}e_{|\Gamma}^{0}||_{\infty}\right)\\
&=&||P^k||_{\infty}
\end{eqnarray}

Hence,
\begin{equation}
\lim_{k->\infty}||P^k||_{\infty}^{\frac{1}{k}} = \rho(P) = \rho(T_{RAS})
\end{equation} $ [] $
\end{proposition}

		\subsection{Convergence of ARAS and ARAS2 in their approximated form} 
		\label{subsection-cv-ARAS-approx}

As we mentioned previously, there exist different approaches to approximate the error transfer operator $P$. For a fully algebraic approach, one will choose the approximation of the operator in a basis built explicitly as we described in \ref{section-approx-explicit}. When it is possible, one can build a complete base analytically and make an approximation of the operator in this base by truncation. In the following, we choose the analytical approach to study the convergence of the method.

Here we focus on elliptic and separable operators. We propose a theorem giving the convergence rate of an ARAS iterative process when the error's transfer operator can be exactly computed in a space spanned by the eigenvectors of $P$ and then truncated to provide an approximation of the error transfer operator in the physical space. 

\begin{theorem}\label{theorem-ARAS-cv}
Let $A$ be a discretized operator of an elliptic problem on a domain $\Omega$. Let us consider a RAS iterative process such as $T_{RAS} = I - M_{RAS}^{-1}A$ defined on $p$ domains.
Let the error transfer operator $P$ on an interface $\Gamma$ be diagonalisable.
If $P$ is diagonalisable, its decomposition in eigenvalues leads to have $P = \mathbb{U} \hat{\hat{P}} \mathbb{U}^{-1}$ where for $i \in \llbracket 1,n \rrbracket$, $\hat{\hat{P}} = diag(\lambda_i)$.
The error on the interface $\Gamma$ in the approximation space follows $\hat{\hat{e}}_{|\Gamma}^{k+1} = \hat{\hat{P}} \hat{\hat{e}}_{|\Gamma}^{k}$.
Each mode converges linearly and independently from the others following $\hat{\hat{e}}_{|\Gamma,i}^{k+1} = \lambda_{i}\hat{\hat{e}}_{|\Gamma,i}^{k}$.
Let $Q_{\lambda} \in \mathbb{R}^{n \times n}$ be a diagonal matrix such that $q_{l} = 1$ if $ 1\leq l \leq q$ and $q_{l} = 0$ if $ q < l $.
And let $\bar{Q}_{\lambda} = I_n -  Q_{\lambda}$.
A coarse approximation of $\hat{\hat{P}}$ can be done choosing a set of $q$ strong modes as $\tilde{P} = Q_{\lambda}\hat{\hat{P}}$.
Writing the preconditioner as:
\begin{equation}\label{eq:ARASq-theorem}
M^{-1}_{ARAS(q),\delta} = \left(I_{m}+R_{\Gamma}^{T}\mathbb{U}\left( \left(I_n-\tilde{P}\right)^{-1}-I_n\right)\mathbb{U}^{-1} R_{\Gamma}\right)\displaystyle\sum_{i=1}^{p}\tilde{R}_{i,\delta}^{T}A_{i,\delta}^{-1}R_{i,\delta} 
\end{equation}
The spectral radius of $T_{ARAS(q)}$ is :
\begin{equation}\label{eq:ARASq-cv}
\rho(T_{ARAS(q)}) = \rho(\bar{Q}_{\lambda} \hat{\hat{P}}) = \lambda_{q+1} < \min \lbrace |\lambda| : \lambda \in \lambda(Q_{\lambda} \hat{\hat{P}}) \rbrace
\end{equation}
\noindent{\bf Proof} We consider the assumptions of the theorem and the formula of the approximated preconditioner \ref{eq:ARASq-theorem}.
The equation \eqref{eq:ARAS} can be written for ARAS(q) such as: 
\begin{eqnarray*}
u^{*} &=& T_{RAS} u^{k-1} + M_{RAS}^{-1}b \\
& &+ R_{\Gamma}^{T}\left(I_n-P\right)^{-1}\left( u_{|\Gamma}^k - Pu_{|\Gamma}^{k-1}\right)R_{\Gamma} \\
& &- R_{\Gamma}^{T} I_n  R_{\Gamma} \left(T_{RAS} u^{k-1} + M_{RAS}^{-1}b \right)
\end{eqnarray*}

We consider that on the interface: 
\begin{equation}\label{eq:RAS-on-gamma}
\left(T_{RAS} u^{k-1} + M_{RAS}^{-1}b\right)_{|\Gamma} = P u_{|\Gamma}^{k-1} + c
\end{equation}
With $c \in \mathbb{R}^n$, a constant vector independent of $u_{|\Gamma}$.

Extracting the interface's solution in the approximation space,
\begin{equation*}
\hat{\hat{u}}_{|\Gamma}^{*} = \hat{\hat{P}} \hat{\hat{u}}_{|\Gamma}^{k-1} + \hat{\hat{c}} + \left(I_n-\tilde{P}\right)^{-1}\left( \hat{\hat{u}}_{|\Gamma}^k - \tilde{P}\hat{\hat{u}}_{|\Gamma}^{k-1}\right)  - Q_{\lambda}\hat{\hat{P}}\hat{\hat{u}}_{|\Gamma}^{k-1} - Q_{\lambda}\hat{\hat{c}}
\end{equation*}
Then,
\begin{equation}
\hat{\hat{u}}_{|\Gamma}^{k} = \bar{Q}_{\lambda}\hat{\hat{P}}\hat{\hat{u}}_{|\Gamma}^{k-1} + Q_{\lambda}\hat{\hat{u}}_{|\Gamma}^{\infty} + \bar{Q}_{\lambda}\hat{\hat{c}}
\end{equation}

The error on the interface is 
\begin{equation*}
\hat{\hat{u}}_{|\Gamma}^{\infty} - \hat{\hat{u}}_{|\Gamma}^{k} = \hat{\hat{u}}_{|\Gamma}^{\infty} - \bar{Q}_{\lambda}\hat{\hat{P}}\hat{\hat{u}}_{|\Gamma}^{k-1} - Q_{\lambda}\hat{\hat{u}}_{|\Gamma}^{\infty} - \bar{Q}_{\lambda}\hat{\hat{c}}
\end{equation*}

Thus,
\begin{equation}
\hat{\hat{e}}_{|\Gamma}^{k} = \bar{Q}_{\lambda}\left( \hat{\hat{u}}_{|\Gamma}^{\infty} - \hat{\hat{P}}\hat{\hat{u}}_{|\Gamma}^{k-1} - \hat{\hat{c}}\right)
\end{equation}

Regarding equation \eqref{eq:RAS-on-gamma}, $\hat{\hat{P}}\hat{\hat{u}}_{|\Gamma}^{k-1} + \hat{\hat{c}}=\hat{\hat{u}}_{|\Gamma}^{k} $, and then,
\begin{equation}
\hat{\hat{u}}_{|\Gamma}^{\infty} - \hat{\hat{P}}\hat{\hat{u}}_{|\Gamma}^{k-1} - \hat{\hat{c}} = \hat{\hat{e}}_{\Gamma}^k = \hat{\hat{P}}\hat{\hat{e}}_{\Gamma}^k
\end{equation}

Hence we write,
\begin{equation}
\hat{\hat{e}}_{|\Gamma}^{k} =  \bar{Q}_{\lambda}\hat{\hat{P}}\hat{\hat{e}}_{|\Gamma}^{k-1}
\end{equation}

We showed that the ARAS iterative process has an error's transfer operator,  $\bar{Q}_{\lambda}\hat{\hat{P}}$, equal to the part of the error's transfer operator $\hat{\hat{P}}$ that we did not compute. We note that $||\bar{Q}_{\lambda}\hat{\hat{P}}|| \leq ||\hat{\hat{P}}|| < 1 $.

$A$ is a discretized operator of an elliptic problem, then we can apply Proposition \ref{prop:convergence-RAS} and write:
\begin{equation}
\rho(T_{ARAS(q)}) = \rho( \bar{Q}_{\lambda}\hat{\hat{P}})
\end{equation}

Then we proof equation \ref{eq:ARASq-cv}. $ [] $
\end{theorem}

\begin{remark}
For a separable operator, the error transfer operator on an interface between two domains is diagonalisable \cite{Garbey1}. Then, the error transfer operator $P$ on an interface $\Gamma$ is diagonalisable for a global operator which is separable and for which the interfaces of all the domains are parallel to one discretization direction. 
\end{remark}

		\subsection{Convergence study in the case of the $2D$  Poisson's equation}		
		\label{subsection-cv-study-2D}
We consider a simple case of an elliptic problem on a rectangle defined in equation \eqref{eq:poisson-2D}. The problem is discretized by $2D$ finite differences and decomposed into $2$ domains $\Omega_1$ and $\Omega_2$ of the same size.

\begin{eqnarray}\label{eq:poisson-2D}
\left\{ \begin{array}{lcl}
- \triangle u &=& f,\; \textrm{in} \;\Omega = [0,1]\times[0,\pi] \\
u&=&0, \; \textrm{on} \;  \partial \Omega \end{array} \right. 
\end{eqnarray}

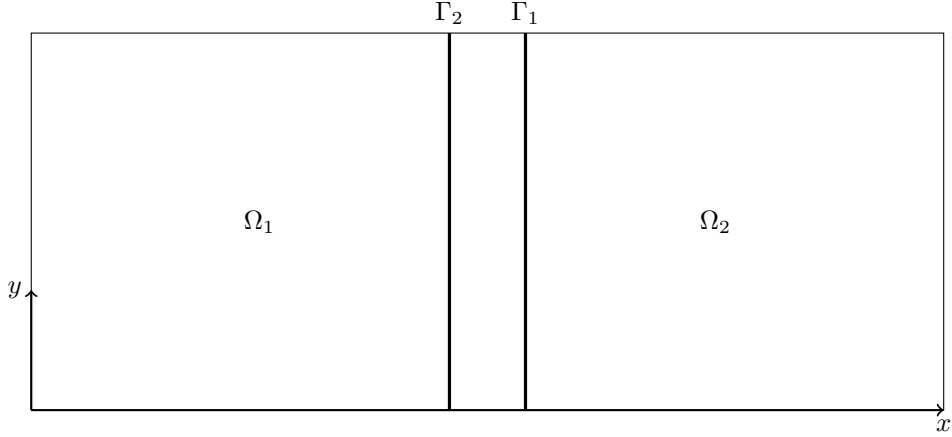
\begin{figure}[ht]
\begin{center}
\input{./figures/figure-2D-domain-poisson.tex}
\caption{$2D$ domain decomposition for Poisson's equation}
\end{center}
\end{figure}

In the Fourier base, for a separable operator with $2$ artificial interfaces, the acceleration is written with $\hat{P} = \bigl( \begin{smallmatrix} 0 & \hat{P}_{\Gamma_2}\\ \hat{P}_{\Gamma_1}&0 \end{smallmatrix}\bigr)$ such as 
\begin{equation}
\begin{pmatrix} \hat{e}_{\Gamma_{1}}^{k}\\ \hat{e}_{\Gamma_{2}}^{k} \end{pmatrix}=\begin{pmatrix} 0 & \hat{P}_{\Gamma_2}\\ \hat{P}_{\Gamma_1} & 0 \end{pmatrix} \begin{pmatrix} \hat{e}_{\Gamma_1}^{k-1}\\ \hat{e}_{\Gamma_2}^{k-1} \end{pmatrix}
\end{equation}

We study the case of a regular grid and write the semi-discretized $2D$ Poisson's operator on $[0,\Gamma_1]\times[0,\pi] \cup [\Gamma_2,1]\times[0,\pi]$, $\Gamma_1 > \Gamma_2$. We consider $N_x+1$ in the direction $x$. The overlap in terms of number of step size in direction $x$ is denoted by $\delta$. And we denote by $\lambda_l$ the eigenvalues of the discretized operator $-\frac{\partial^2}{\partial y^2}$ considering a second order finite difference scheme. We find the coefficient of the matrices $\hat{P}_{\Gamma_i}$ solving the equations:
\begin{eqnarray}
\left\{ \begin{array}{c c c} - \frac{\hat{u}_{i+1,l}^1-2\hat{u}_{i,l}^1 + \hat{u}_{i-1,l}^1}{h_{x}^2} + \lambda_l\hat{u}_{i,l}^1 &=& 0 \text{ , } i \in \llbracket 1, N_x-1 \rrbracket \\ 
\hat{u}_{0,l}^1&=&0 \\
\hat{u}_{N_x,l}^1&=&1 \end{array} \right. \\
\left\{ \begin{array}{c c c} - \frac{\hat{u}_{i+1,l}^2-2\hat{u}_{i,l}^2 + \hat{u}_{i-1,l}^2}{h_{x}^2} + \lambda_l\hat{u}_{i,l}^2 &=& 0 \text{ , } i \in \llbracket 1, N_x-1 \rrbracket \\ 
\hat{u}_{0,l}^2&=&1 \\
\hat{u}_{N_x,l}^2&=&0  \end{array} \right.
\end{eqnarray}

The roots of those equations are such as,
\begin{eqnarray}
r_1 &=& \frac{2 + \lambda_l h_{x}^2 + \sqrt{\lambda_{l}^2 h_{x}^4 + 4 \lambda_l h_{x}^2}}{2}\\
r_2 &=& \frac{2 + \lambda_l h_{x}^2 - \sqrt{\lambda_{l}^2 h_{x}^4 + 4 \lambda_l h_{x}^2}}{2}
\end{eqnarray}

Then for the first domain $\Omega_1$ the solutions have the form:
\begin{equation}
\hat{u}_{l,j}^1 = \frac{r^{j}_1-r^{j}_2}{r^{N_x}_1-r^{N_x}_2}
\end{equation}

Thus,
\begin{equation}
\hat{u}_{N_x-\delta,l}^1 = \frac{r^{N_x-\delta}_1-r^{N_x-\delta}_2}{r^{N_x}_1-r^{N_x}_2}
\end{equation}

$\hat{P}_{\Gamma_1}$ appears to be diagonal and its diagonal coefficients $\delta_{l,\Gamma_1}$ can be analytically derived such as: 
\begin{eqnarray}
\delta_{l,\Gamma_1} &=& \hat{u}_{N_x-\delta,l}^1* \hat{u}_{\delta,l}^2\\
& = & \left( \frac{r^{N_x-\delta}_1-r^{N_x-\delta}_2}{r^{N_x}_1-r^{N_x}_2} \right)\left( \frac{-r^{N_x}_2r^{\delta}_1+r^{N_x}_1r^{\delta}_2}{r^{N_x}_1-r^{N_x}_2} \right)
\end{eqnarray}

Exactly the same development can be done to find $\delta_{l,\Gamma_2}$. 

\begin{remark}
If the two domains have the same size, then $\delta_{l,\Gamma_1}=\delta_{l,\Gamma_2}$.
\end{remark}

The $\hat{P}$ matrix has the form:
\begin{equation}\label{eq:P-2D-analytic}
\hat{P} = \begin{pmatrix}
0                                     &              &     & \delta_{1,\Gamma_2} &             &                                          \\
                                       &  \ddots &     &                                        & \ddots &                                          \\
                                       &               & 0 &                                        &              &  \delta_{n,\Gamma_2}  \\
\delta_{1,\Gamma_1} &               &     & 0                                    &              &                                          \\
                                       & \ddots   &     &                                       & \ddots &                                          \\
                                       &               & \delta_{n,\Gamma_1} &    &              & 0                           
\end{pmatrix}\text{ , with } 1>\delta_{1,\Gamma_i}\geq ... \geq \delta_{n,\Gamma_i}.
\end{equation}

For the sake of simplicity we consider that the domain $\Omega_1$ and $\Omega_2$ have the same size. We note $\delta_{l} = \delta_{l,\Gamma_1} = \delta_{l,\Gamma_2}$.
Then we can calculate the determinant of $(\hat{P}-\lambda I_n)$: 
\begin{eqnarray*}
det(\hat{P}-\lambda I_n) &=& \prod_{l=1}^{n}(\delta_{l}-\lambda)(\delta_{l}+\lambda)
\end{eqnarray*}
Hence the spectrum is $\left\{\delta_1,...,\delta_n,-\delta_1,...,-\delta_n \right\}$.

We showed that the error's transfer operator is diagonalisable for this problem. $P$ can be written in a base $\mathbb{U}$ of eigenvectors as follows: 
\begin{equation}
P = \mathbb{U} \hat{\hat{P}} \mathbb{U}^{-1} \text{, with }  \hat{\hat{P}} = diag(\delta_1,...,\delta_n,-\delta_1,...,-\delta_n )
\end{equation}

Then we can estimate the convergence rate of the RAS, ARAS(q) and ARAS2(q) applying Theorem \ref{theorem-ARAS-cv}.
\begin{eqnarray}\label{eq:theoretical-results-cv-2D}
\rho(T_{RAS}) &=& \delta_1 \\
\rho(T_{ARAS(q)}) &=& \delta_{q+1} \\
\rho(T_{ARAS2(q)}) &=& \delta_{q+1}^2
\end{eqnarray}

Because the eigenvalues and the values of $\hat{P}_{\Gamma_i}$ are equal we can verify a correspondence between the approximation by truncation in the eigenvectors space and the Fourier space. Selecting the first Fourier mode corresponds to selecting the highest eigenvalues.
Let us introduce the transfer matrix $C_{\Gamma_i}$ from the real space to the Fourier space and the transfer matrix $D_{\Gamma_i}$ from the Fourier space to the real space:

\begin{equation*}
\begin{array}{c c c  c c c c}
C_{\Gamma_i} & : & \mathbb{R}^n \longrightarrow \mathbb{C}^n & \text{ and } & D_{\Gamma_i} & : & \mathbb{C}^n \longrightarrow \mathbb{R}^n\\
& & e_{|\Gamma_i} \longmapsto \hat{e}_{|\Gamma_i}& & & & \hat{e}_{|\Gamma_i} \longmapsto e_{|\Gamma_i}
\end{array}
\end{equation*}

Then we write
\begin{equation}\label{eq:P-in-Fourier-base}
P = D\hat{P}C = \begin{pmatrix} D_{\Gamma_1} & 0 \\ 0 & D_{\Gamma_2} \end{pmatrix} \begin{pmatrix} 0 & \hat{P}_{\Gamma_2} \\ \hat{P}_{\Gamma_1} & 0 \end{pmatrix} \begin{pmatrix} C_{\Gamma_1} & 0 \\ 0 & C_{\Gamma_2} \end{pmatrix}
\end{equation}
The approximation is done by applying the operator $$Q_{\cal{F}} = \begin{pmatrix}Q_{\Gamma_{1},\cal{F}} & 0\\0 &  Q_{\Gamma_{2},\cal{F}} \end{pmatrix} $$ where $Q_{\Gamma_{1},\cal{F}}=diag(q_{l})$, $q_{l} = 1$ if $ 1\leq l \leq q$ and $q_{l} = 0$ if $ q < l $.
Then we write the preconditioner 
\begin{equation}
M_{ARAS(q)}^{-1} = (I + R_{\Gamma}^{T}D(I_n-Q_{\cal{F}}\hat{P})^{-1}-I_n)CR_{\Gamma})M_{RAS}^{-1}
\end{equation}

As previously we introduce a matrix $ \bar{Q}_{\cal{F}} = (I - Q_{\cal{F}}) $.

We can then follow the demonstration done in the proof of theorem \ref{theorem-ARAS-cv} writing the error on the interface in the Fourier space as 
\begin{equation}
\hat{u}_{|\Gamma}^{\infty} - \hat{u}_{|\Gamma}^{k} = \hat{u}_{|\Gamma}^{\infty} - \bar{Q}_{\cal{F}}\hat{P}\hat{u}_{|\Gamma}^{k-1} - Q_{\cal{F}}\hat{u}_{|\Gamma}^{\infty} - \bar{Q}_{\cal{F}}\hat{c}
\end{equation}

Thus,
\begin{equation}
\hat{e}_{|\Gamma}^{k} = \bar{Q}_{\cal{F}}\left( \hat{u}_{|\Gamma}^{\infty} - \hat{P}\hat{u}_{|\Gamma}^{k-1} - \hat{c}\right)
\end{equation}

Regarding equation \eqref{eq:RAS-on-gamma}, $\hat{P}\hat{u}_{|\Gamma}^{k-1} + \hat{c}=\hat{u}_{|\Gamma}^{k} $, and then,
\begin{equation}
\hat{e}_{|\Gamma}^{k} =  \bar{Q}_{\cal{F}}\hat{e}_{|\Gamma}^{k} 
\end{equation}

As $\hat{e}_{|\Gamma}^{k} = \hat{P}\hat{e}_{|\Gamma}^{k-1}$ we write
\begin{equation}
\hat{e}_{|\Gamma}^{k} =  \bar{Q}_{\cal{F}}\hat{P}\hat{e}_{|\Gamma}^{k-1}
\end{equation}

We showed that the ARAS iterative process has an error transfer operator,  $\bar{Q}_{\cal{F}}\hat{P}$, equal to the part of the error transfer operator $\hat{P}$ that we did not compute. We note that $||\bar{Q}_{\lambda}\hat{P}|| \leq ||\hat{P}|| < 1 $.

We can apply Proposition \ref{prop:convergence-RAS} and write:
\begin{equation}
\rho(T_{ARAS(q)}) = \rho( \bar{Q}_{\cal{F}}\hat{P})
\end{equation}

The conclusion becomes the same as applying Theorem \ref{theorem-ARAS-cv}.

We pointed out here the way the approximation of the error's transfer operator affects the convergence of Schwarz iterative processes in the case of a separable operator for a two domain decomposition. This enables us to understand the philosophy of approximating the matrix $P$  in different spaces and links the works done in \cite{Garbey1,Baranger,mg_dtd_fd}.

\section{Results on academic problems}
\label{section-results-academic}

		\subsection{$2D$ theoretical study}
		\label{subsection-numeric-2D-theo}

The goal of this section is to validate the ARAS method on a simple case where the ARAS preconditioner can be written analytically and where we can apply Theorem \ref{theorem-ARAS-cv}. 
We consider the $2D$ problem decomposed in $2$ domains presented in subsection \ref{subsection-cv-study-2D}. The grid size is about $32 \times 32$. This subsection provides the theoretical framework we implement in Matlab. Here, we verify numerically the theoretical results given previously.

\begin{figure}[ht]
\begin{center}
\includegraphics[width=100mm,keepaspectratio]{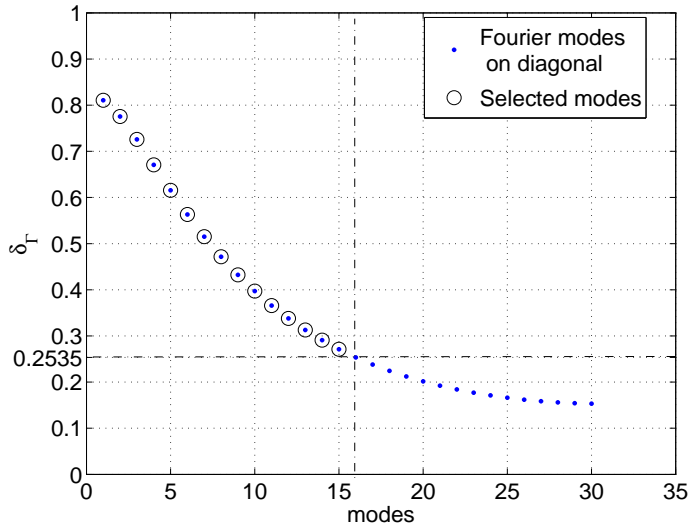}
\caption{Diagonal coefficients of $\hat{P}$ and $Q_{\cal{F}}\hat{P}$ corresponding to the modes to accelerate.}\label{figure-numeric-2D-fourier-modes}
\end{center}
\end{figure}

We build the matrix $\hat{P} \in \mathbb{C}^{30 \times 30}$. Only the internal points are taken, leading to $30$ modes which can be accelerated. Those modes decrease from $0.8106$ to $0.1531$. We decide to compute the entire $\hat{P}$ and a truncated one of size $q=15$, $Q_{\cal{F}}\hat{P}$. Figure \ref{figure-numeric-2D-fourier-modes} shows the coefficient of the matrix computed.
The goal here is to retrieve the convergence rate given by the application of theorem \ref{theorem-ARAS-cv} in equation \eqref{eq:theoretical-results-cv-2D}. The convergence rate of a ARAS(q) type preconditioner is related to the coefficients of $\hat{P}$ denoted by $\delta_l$.

For $q=15$ we note that, 
\begin{eqnarray}
\delta_1 &=&  0.8106\\
\delta_{q+1} &=& 0.2535 \\
\delta_{q+1}^2 &=& 0.0643
\end{eqnarray}

Then we compute the RAS, ARAS(q) and ARAS2(q) preconditioners and compute for each preconditioner of type $*$ the spectral radius of $T_*$ and the conditioning number in the $2$ norm of $M_{*}^{-1}A$. The convergence for each appropriate stopping criteria is $10^{-10}$. 

\begin{table}[h]
\begin{center}
\begin{tabular}{|c|c|c|c|c|}
\hline
prec. $*$ & $\rho(T_*)$ & $ \kappa(M_{*}^{-1}A) $ & It. Rich. & It. GCR \\
\hline
RAS                       & 0.8106          & 30.0083                       & 96            & 18 \\
ARAS(q=15)        &  0.2535         & 5.2358                          & 14            & 7 \\
ARAS2(q=15)      &  0.0643         & 1.1451                          & 7               & 5 	\\
ARAS2(q=30)      &  1.4319 e-13 & 1.0000                         & 1               & 1 	\\
\hline
\end{tabular}
\caption{Numerical performance of RAS, ARAS and ARAS2 on the $2D$ Poisson problem.}
\label{table-numeric-2D-rho-cond}
\end{center}
\end{table}

Table \ref{table-numeric-2D-rho-cond} shows the numerical convergence rate of RAS, ARAS(15), ARAS2(15), and ARAS2(30) for the test case presented in Subsection \ref{subsection-cv-study-2D}. The numerical values obtained for $\rho{T_*}$ match perfectly the theoretical results. It exhibits also that the Aitken acceleration of the RAS enhances the condition number of the preconditioned problem. In accordance with the theory, when $q=30$, the size of the artificial interface, $P$ is exact and $M_{ARAS2}^{-1} = A^{-1}$ numerically. 

Moreover the number of iterations of the iterative process ARAS(q=15) is twice the number of iterations of the iterative process ARAS2(q=15).

		\subsection{Observing the influence of the partitioning and the approximation space on a $2D$ Helmholtz problem} 
		\label{subsection-numeric-2D-helmohltz}
		
		Here, we focus on the influence of the partitioning chosen to set up the domain decomposition method. We also focus on the influence of the choice of a base to approximate the Aitken's acceleration. 
When the mesh is known it is possible to partition the operator following a geometric partitioning. One point is to see what can happen if we partition the operator with a graph partitioning approach such as METIS. Another point is to see how the choice of a base influence the performance and the cost of the ARAS type preconditioner.
		
Let us consider the 2D Helmholtz problem  $(-\omega  -\triangle)u = f \mbox{ in } \Omega=[0,1]^2, \, u=0 \,\textrm{on}\, \partial \Omega$. 
The problem is discretized by second order finite differences with $m$ points in each direction $x$ and $y$ giving a space step $h=\frac{1}{m-1}$.
The set value $\omega= 0.98 \frac{4}{h^2} (1 - cos(\pi h))$ is close to the minimum eigenvalue of the discrete $-\triangle$ operator in order to have an ill-conditioned discrete problem with $\kappa_\infty(A)=1.7918\, E+07$  for $m=164$. 

\begin{figure}[ht]
\begin{center}
\begin{minipage}{13cm}
\includegraphics[width=60mm,keepaspectratio]{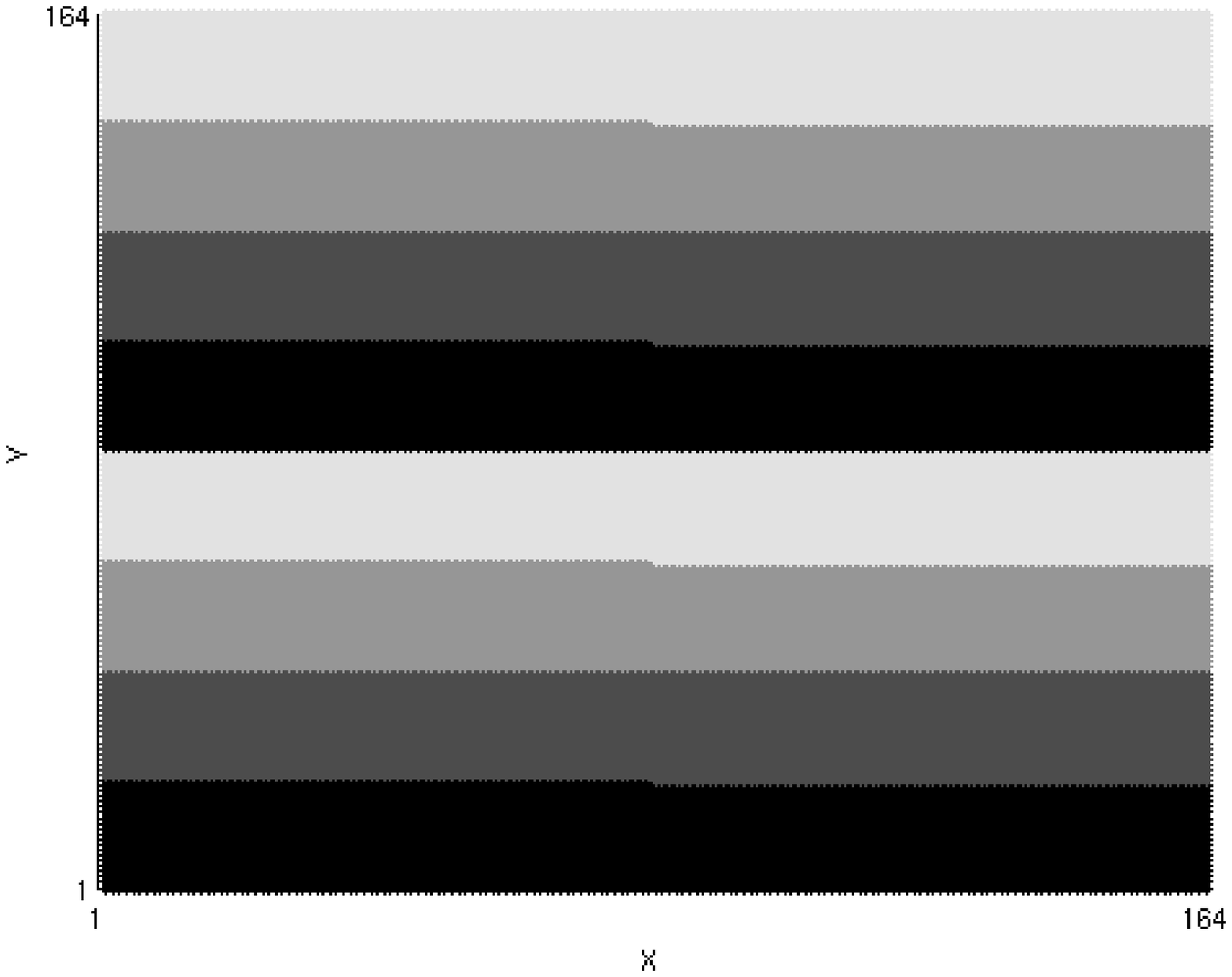}
\includegraphics[width=60mm,keepaspectratio]{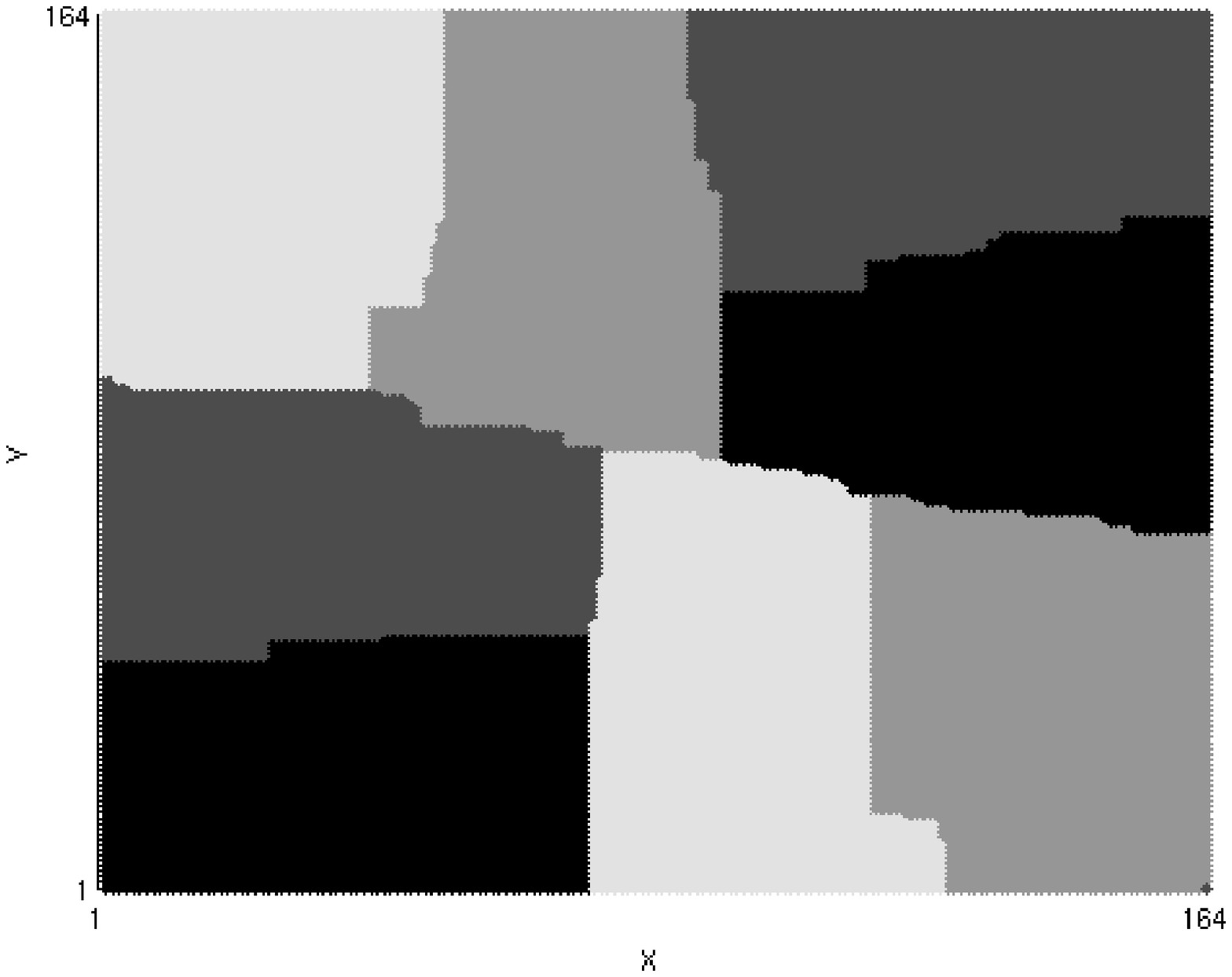}
\end{minipage}
\caption{Partitioning into 8 parts on a 2D Helmholtz problem of size $164 \times 164$, (left) Physical Band partitioning, (right) METIS partitioning.}
\label{partitioning}
\end{center}
\end{figure}

First we solve the problem with a physical band partitioning, and then we solve it with a METIS partitioning using eight sub-domains.
Figure \ref{partitioning} illustrates the physical band and the METIS partitioning using eight sub-domains.
In the physical partitioning, borders are smooth, contrary to the METIS partitioning which creates corners and irregular borders.
The corners give cross points which deteriorate the convergence of the Schwarz method.

For each partitioning, we build the ARAS2 preconditioner in two different bases:
\begin{itemize}
\item an orthogonal base arising from the application of the preconditioner RAS on a sequence of random vectors (see subsection \ref{subsection-random-base}).
\item a base built from the successive Schwarz solution on the interface and passed in its SVD base (see subsection \ref{subsection-SVD-base}).
\end{itemize}

\begin{remark}
In the following, we denote by ARAS2(r=$\frac{n}{q}$) the preconditioner approximated in the "random" base. Because the number of vectors can be high for this kind of base, we choose to express the reduction number $r$ in parentheses instead of $q$, the number of column of $\mathbb{U}_q$, but the formula is still:
\begin{equation*}
M^{-1}_{ARAS(r =\frac{n}{q} ),\delta} = \left(I_{m}+R_{\Gamma}^{T}\mathbb{U}_q\left( \left(I_q-\tilde{P}_{\mathbb{U}_q}\right)^{-1}-I_q\right)\mathbb{U}_q^{-1} R_{\Gamma}\right)\displaystyle\sum_{i=1}^{p}\tilde{R}_{i,\delta}^{T}A_{i,\delta}^{-1}R_{i,\delta} 
\end{equation*}
\end{remark}
\begin{figure}[h]
\begin{center}
\begin{minipage}{13cm}
	\includegraphics[height=80mm,width=60mm]{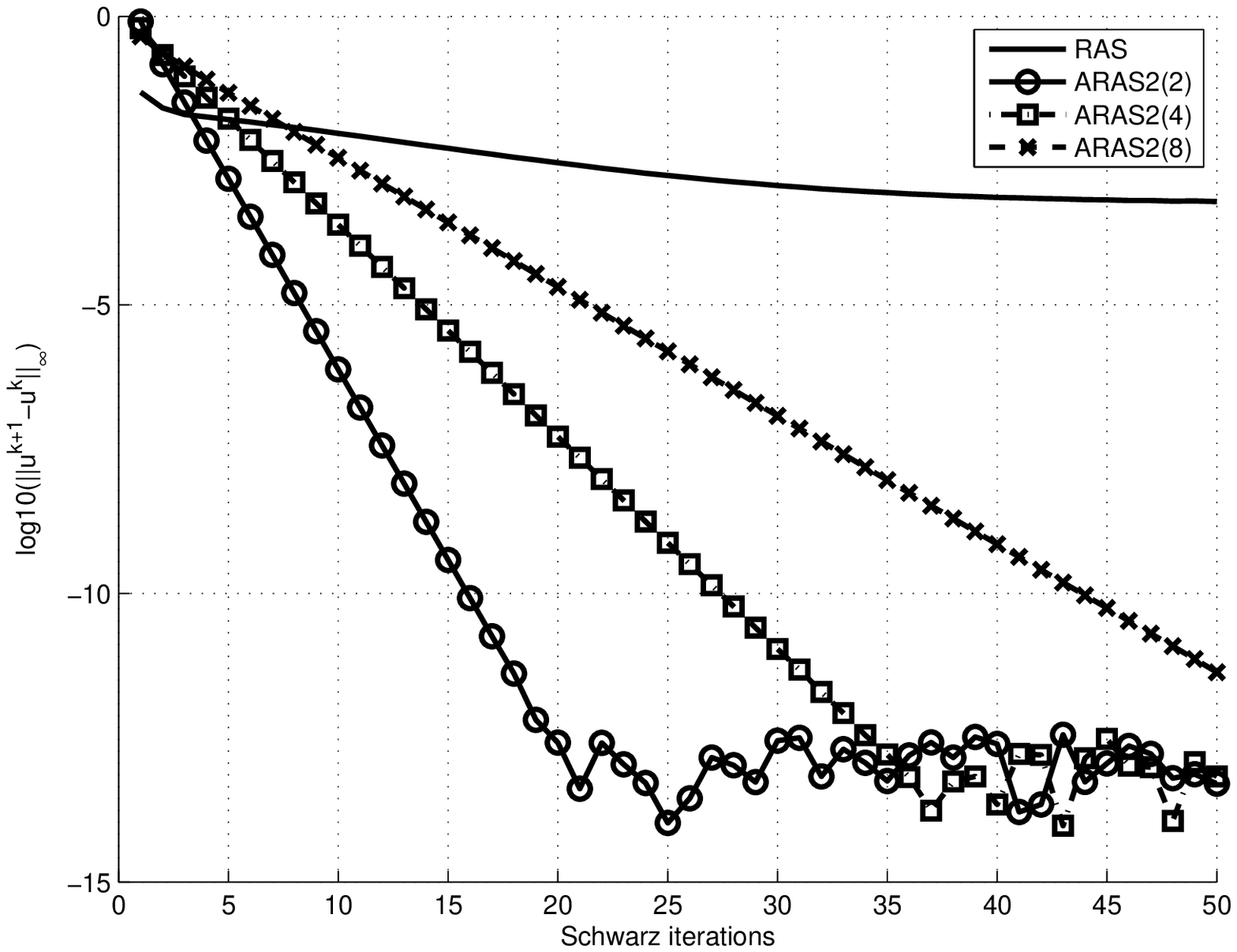}
	\includegraphics[height=80mm,width=60mm]{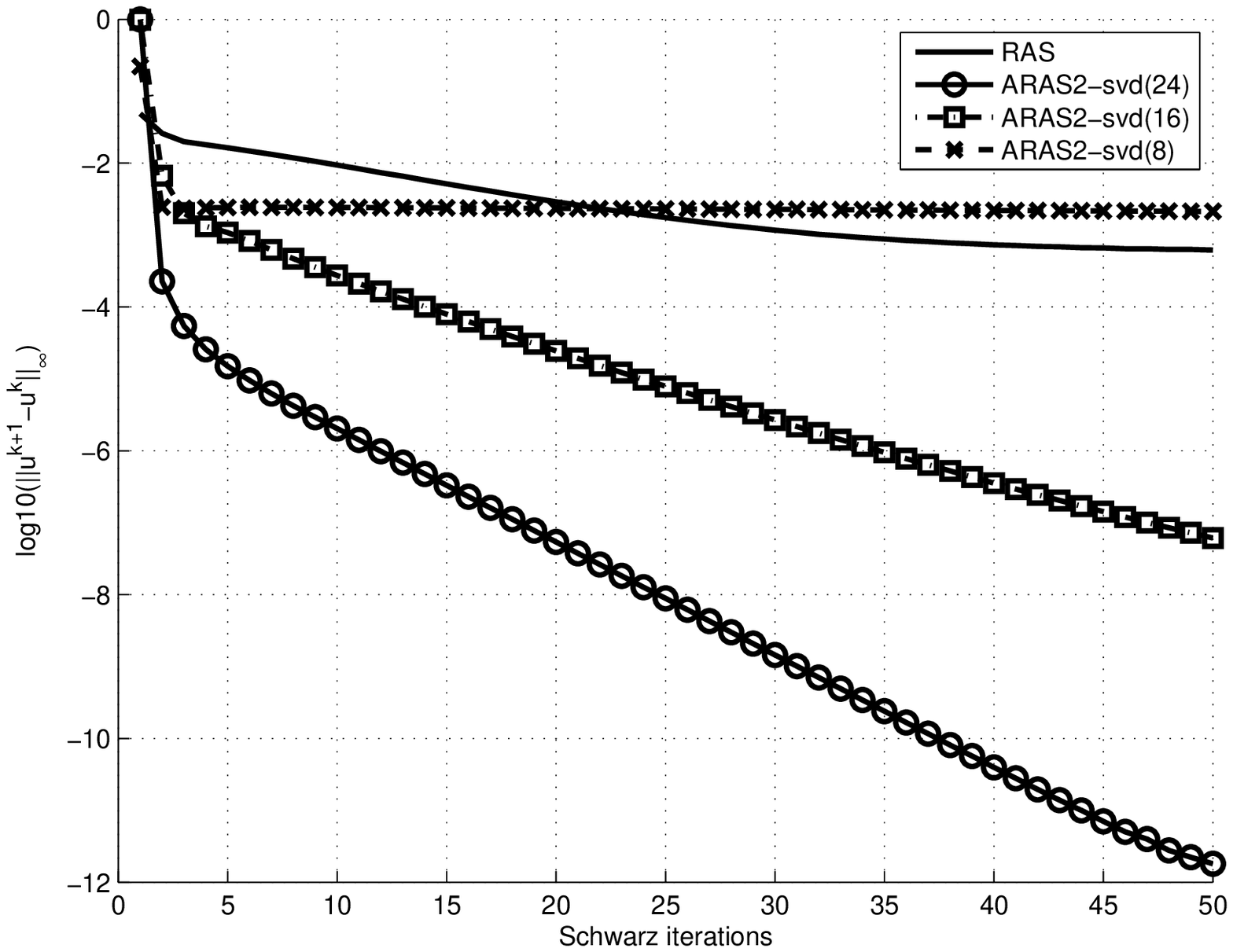}

	\includegraphics[height=80mm,width=60mm]{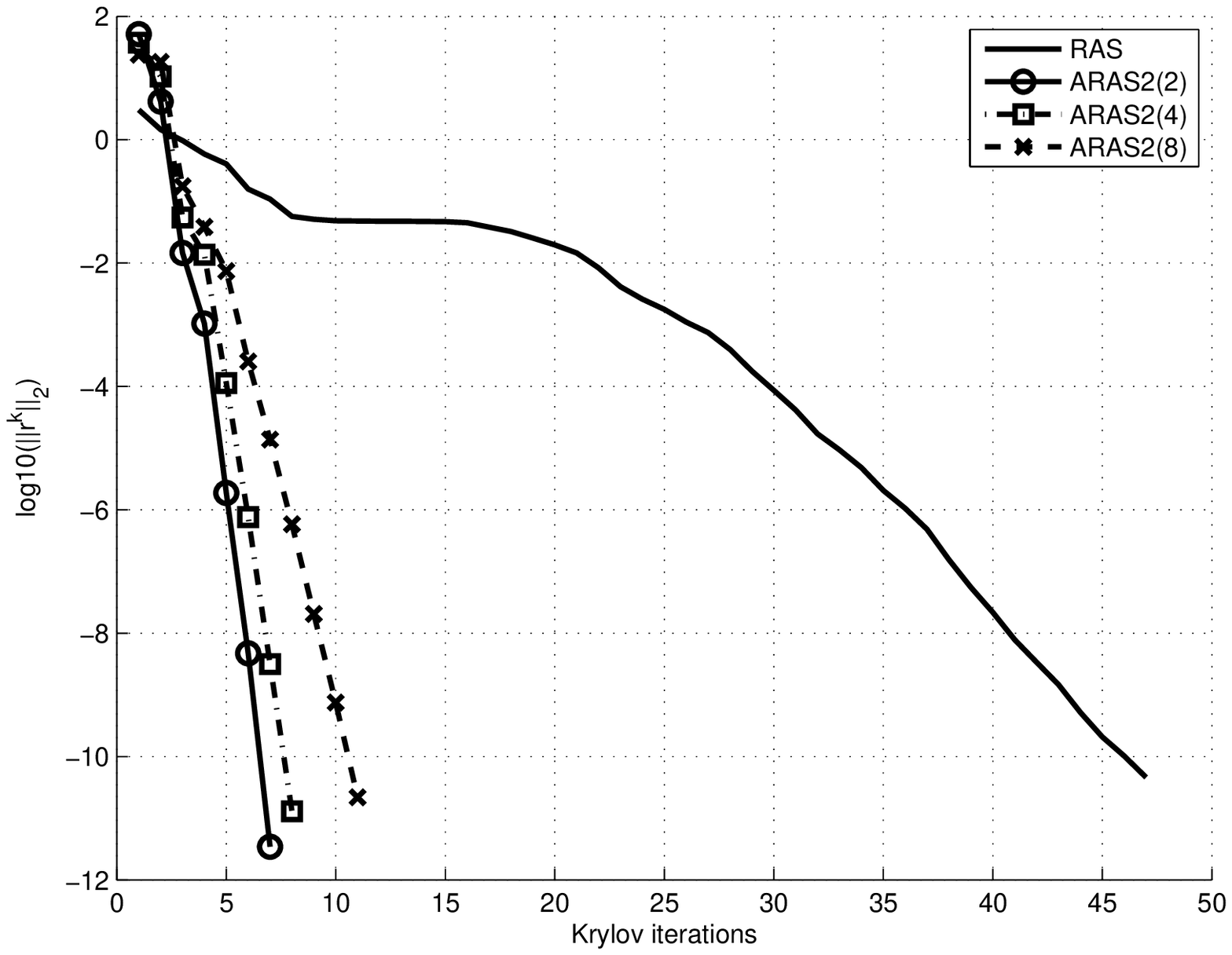}
	\includegraphics[height=80mm,width=60mm]{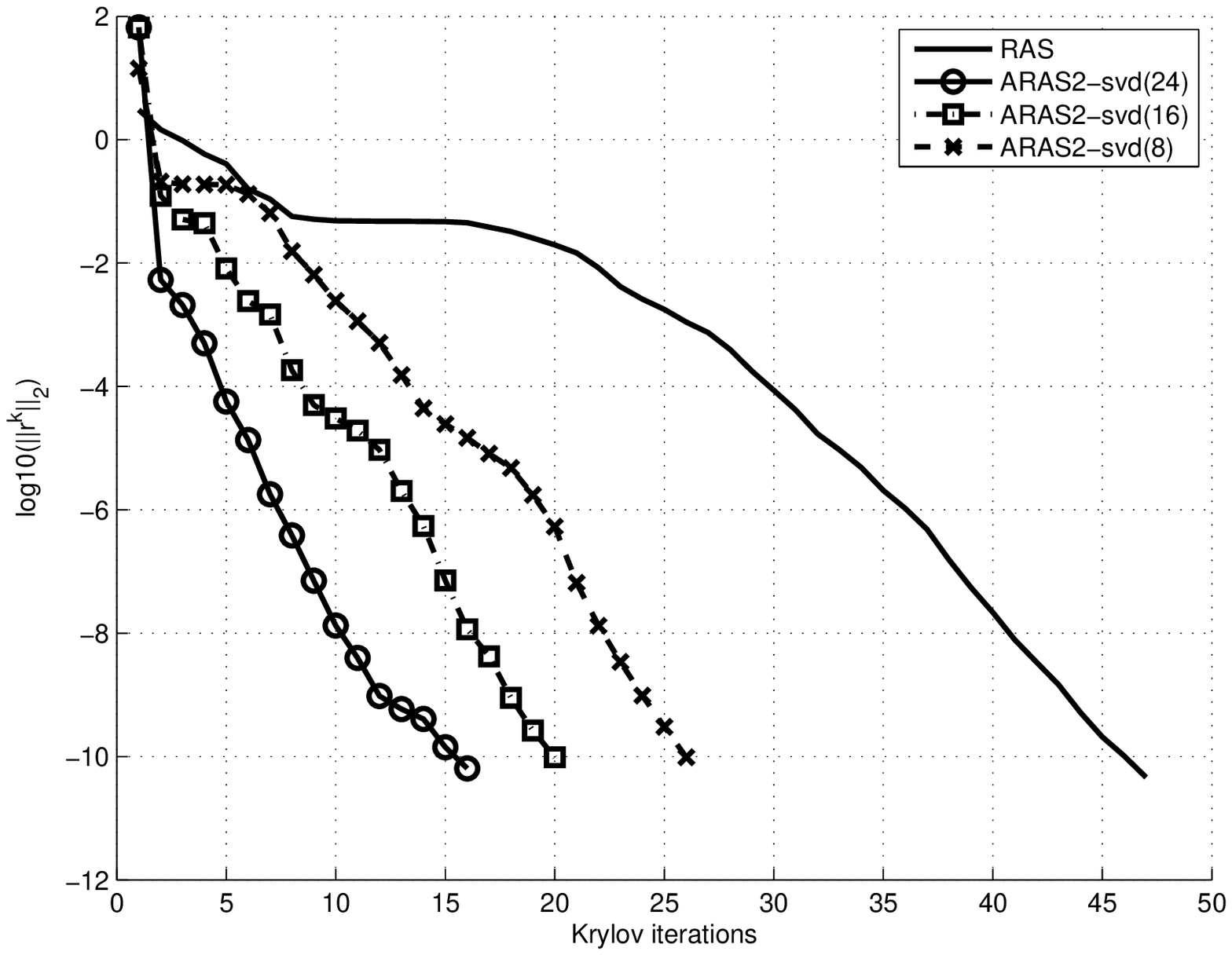}
\end{minipage}
\caption{Solving 2D Helmholtz equation on a $164 \times 164$ Cartesian grid, physical band partitioning, $p =8$,(left) ARAS2(r = $ \frac{n}{q} $ ) is built with a Random base, (right) ARAS2(q) is built with a SVD base, (top) Convergence of Iterative  Schwarz Process, (bottom) convergence of GCR method preconditioned by RAS and ARAS2.}
\label{physicalPartResult}
\end{center}
\end{figure}

\begin{figure}[h]
\begin{center}
\begin{minipage}{13cm}
	\includegraphics[height=80mm,width=60mm]{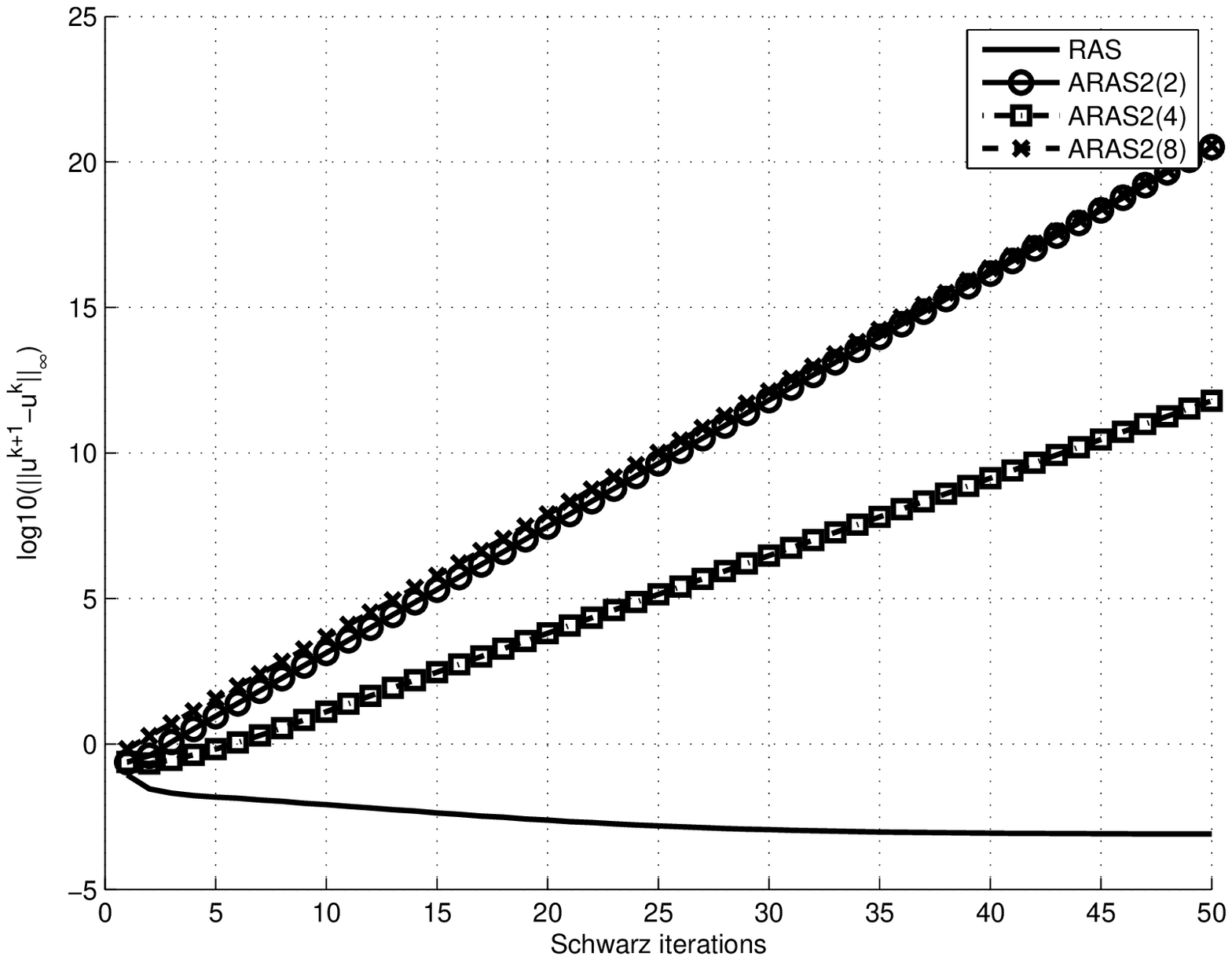}
	\includegraphics[height=80mm,width=60mm]{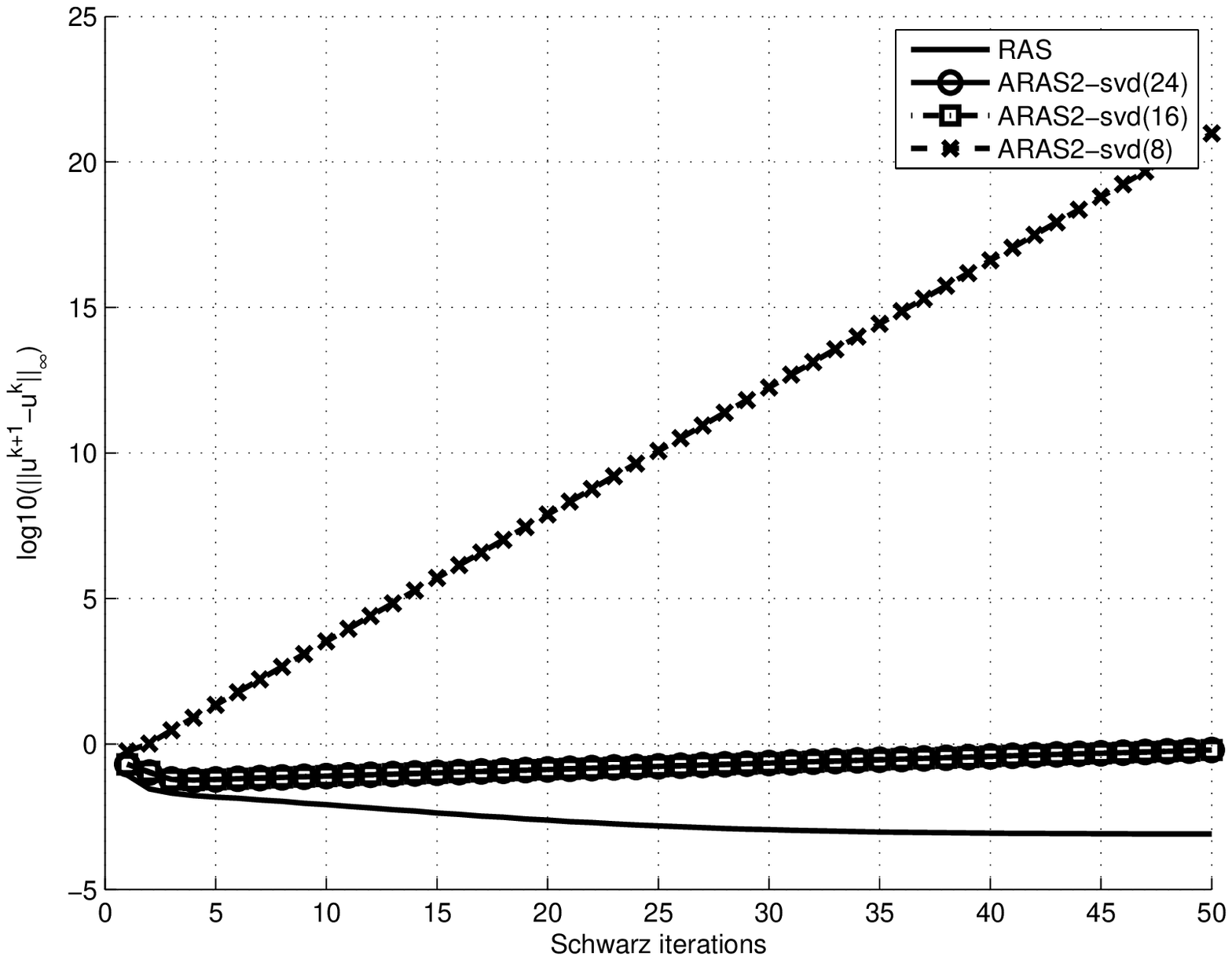}

	\includegraphics[height=80mm,width=60mm]{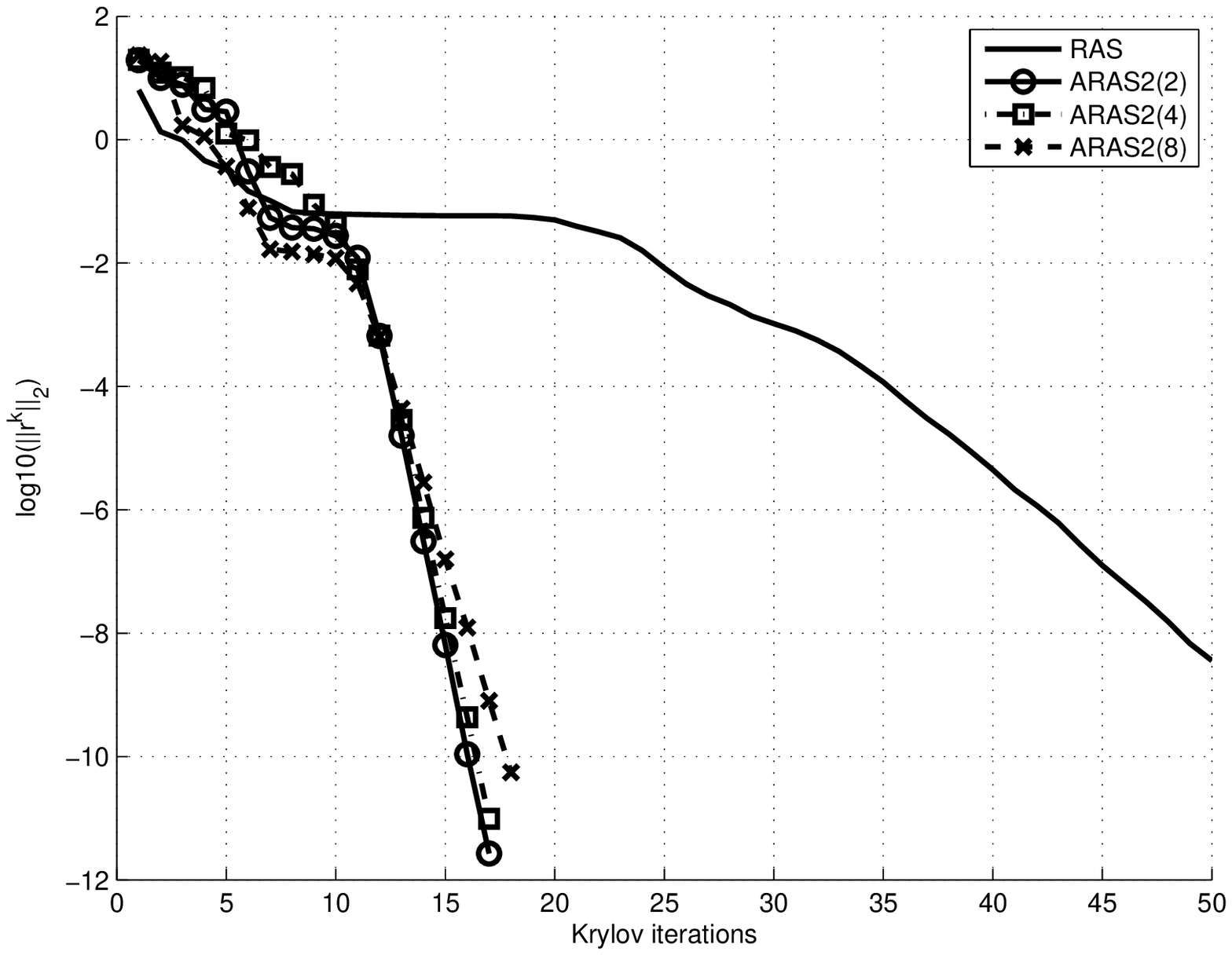}	
	\includegraphics[height=80mm,width=60mm]{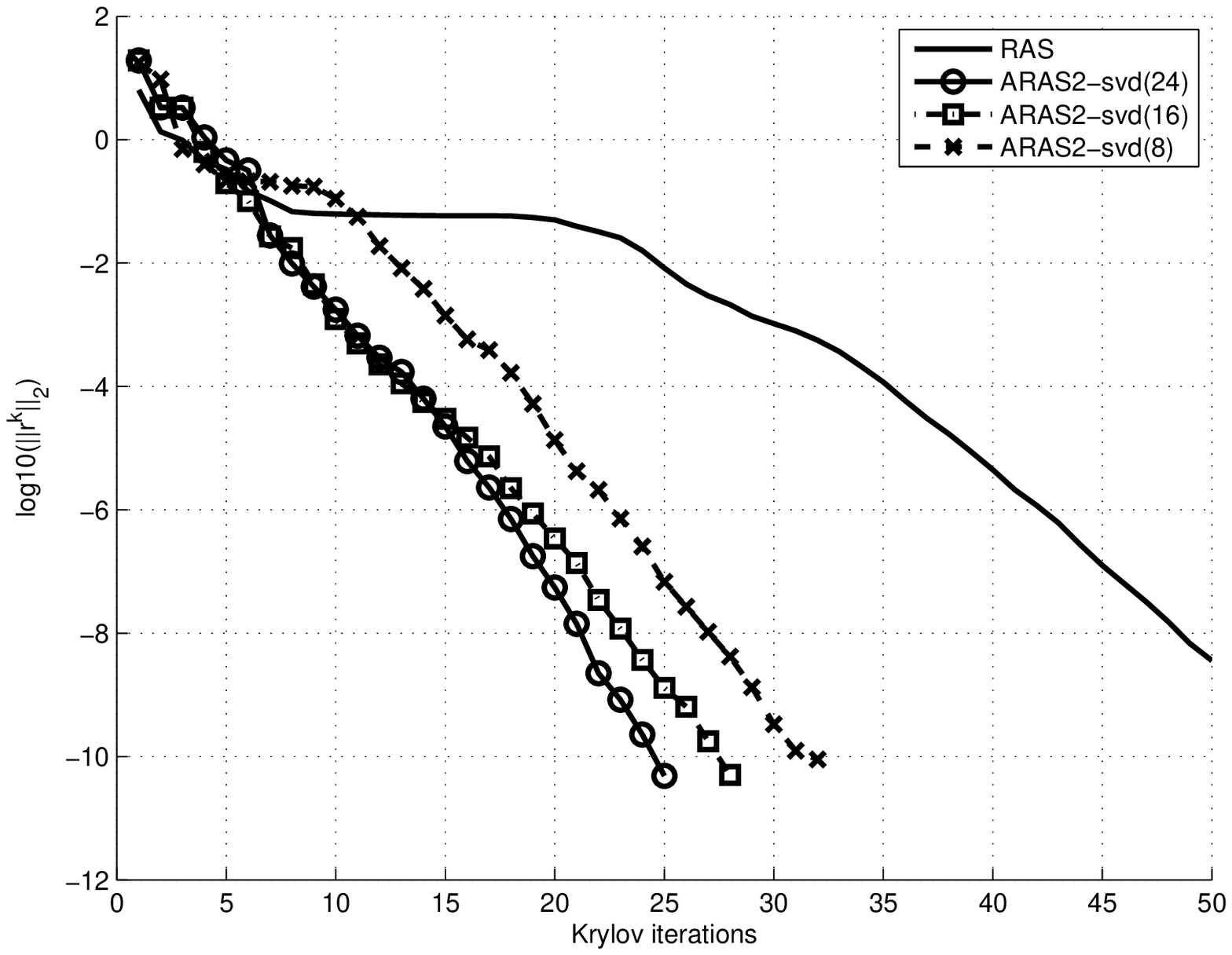}	
\end{minipage}	
\caption{Solving 2D Helmholtz equation on a $164 \times 164$ Cartesian grid, METIS partitioning, $p =8$,(left) ARAS2 is built with a Random base, (right) ARAS2 is built with a SVD base, (top) Convergence of Iterative  Schwarz Process, (bottom) convergence of GCR method preconditioned by RAS and ARAS2.}
\label{metisPartResult}
\end{center}
\end{figure}

Figure \ref{physicalPartResult} (respectively Figure \ref{metisPartResult}) presents the Richardson process with ARAS2 and the ARAS2 preconditioned GCR Krylov method for the physical band partitioning ( respectively a METIS partitioning) using a random base or a SVD base. These results  were obtained with a sequential Matlab code able to run small academic problems. The Krylov method used is the gradient conjugate residual GCR while the LU factorisation is used to solve the sub-domain's problems. 

These results exhibit:
\begin{itemize}
\item Richardson processes in Figure \ref{physicalPartResult} converge while in Figure \ref{metisPartResult} only the RAS iterative process converges and the ARAS2 process diverges. Consequently the physical partitioning enables good convergence of the RAS    which can be accelerated with the approximation of $P$ by $P_{\mathbb{U}_q}$.  The METIS partitioning slows the convergence of the RAS due to the cross points.
Then, using a domain decomposition method as a solver with an algebraic partitioning can produce bad results when it is accelerated by Aitken. Let us notice that the full $P$ makes the RAS process converge in one iteration.
\item Nevertheless, the Aitken-RAS used as a preconditioner is very efficient even on the METIS partitioning with cross points where the Aitken-RAS as a Richardson process diverges. This makes the Aitken-RAS a robust algebraic preconditioner. We must notice that the effect of the preconditioning with a METIS partitioning is less efficient than the one with the physical partitioning.
\item The better the base $\mathbb{U}_q$ is able to represent the interface solution, the better the preconditioner is for the random base and the SVD base.
\end{itemize}

Let us observe the difference between the two choices of base to compute the acceleration. 
On the one hand, the choice of an orthogonal base arising from the application of the RAS preconditioner on a sequence of random vectors presents good advantages for a preconditioner. With this approach, the preconditioner can be used for different right-hand sides. However, the number of vectors necessary to describe the interface can be close to the size of the global interface, increasing the cost of the preconditioner. 
On the other hand, it is possible to build the acceleration for many iterations of the Additive Schwarz process, computing the SVD of the interface solutions. Then  the acceleration process is problem-dependent, but experience shows that a small number of iterations can enable a good approximation of $P_{\mathbb{U}_q}$.

For a physical partitioning, we can evaluate the cost of each preconditioner. For each sub-domain, the artificial interface is of size $164$ for the uppermost or lowermost sub-domain, and $164*2 = 328$ for internal sub-domains. Hence, for $p=8$ the global interface has a size of $164*(6*2+2) = 2296$. For $r=1$ the base is complete and $P_{\mathbb{U}_q}$ is exact. There is no need to use ARAS as a preconditioning technique. 
For all $r$ the size of $\Gamma$ is $\frac{2296}{r}$. Then the number of $M^{-1}_{RAS}x = y$ products to build $M^{-1}_{ARAS}$ is $3*\frac{2296}{r}$. While the number of $M^{-1}_{RAS}x = y$ products  to build $M^{-1}_{ARAS}$ for the base arising from SVD only depends on the number of Richardson iterations.\\
 Figure \ref{physicalPartResult}, shows that for $r=8$, $n=185$ and the number of products $M^{-1}_{RAS}x = y$ is $555$, the convergence of GCR is reached in $11$ iterations. For $24$ iterations of Schwarz, we build a matrix $P_{\mathbb{U}_q}$ of size $24$, which is around eight times smaller than with the previous base.  The number of matrix products is $48$, $12$ times smaller than with the previous base. The number of GCR iterations is $15$. Eventually, the cost of a good independent preconditioner is excessive compared to the one with the SVD.

\begin{figure}[ht]
\begin{center}
\begin{minipage}{9cm}
	\includegraphics[height=80mm,width=90mm]{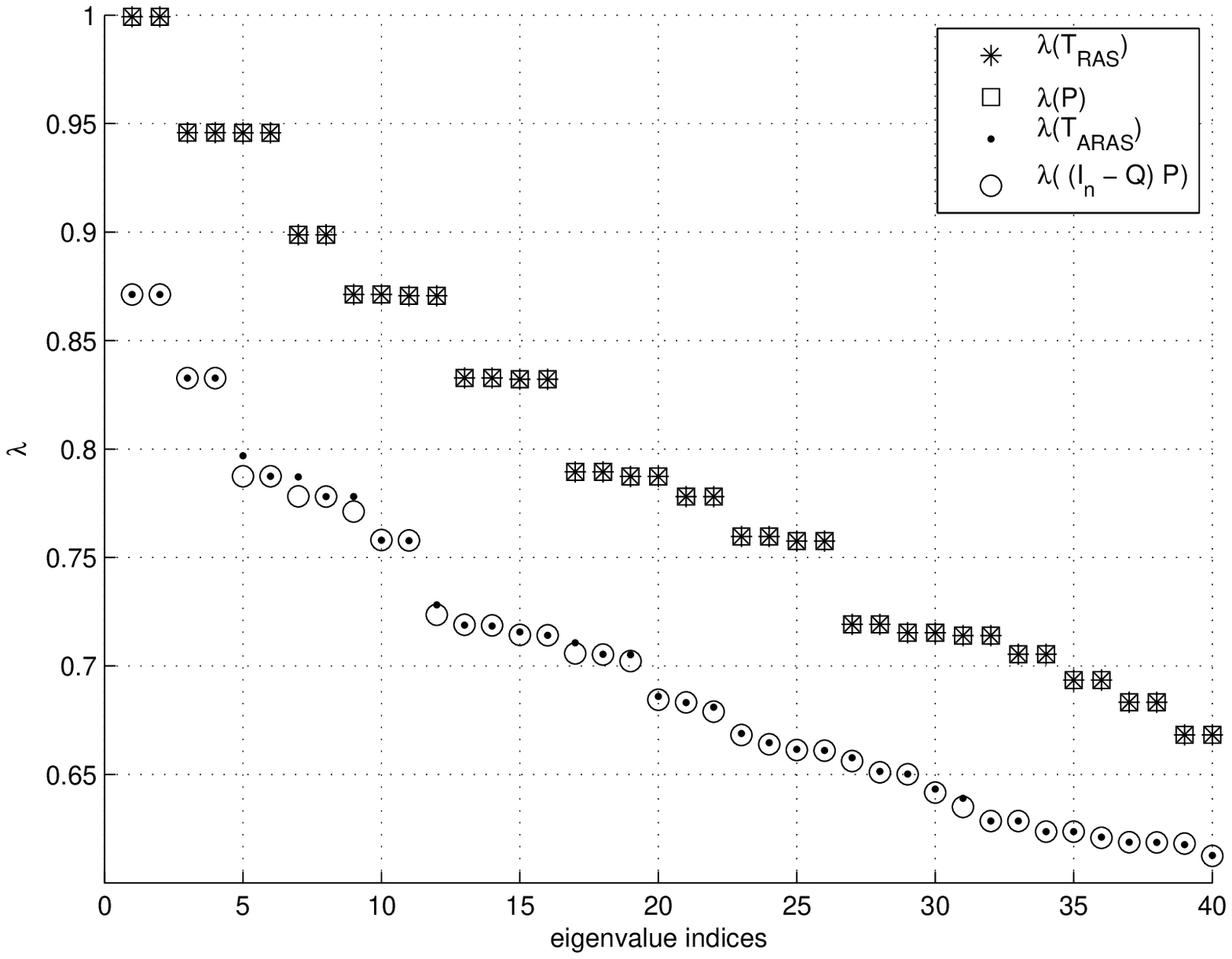}
	\includegraphics[height=80mm,width=90mm]{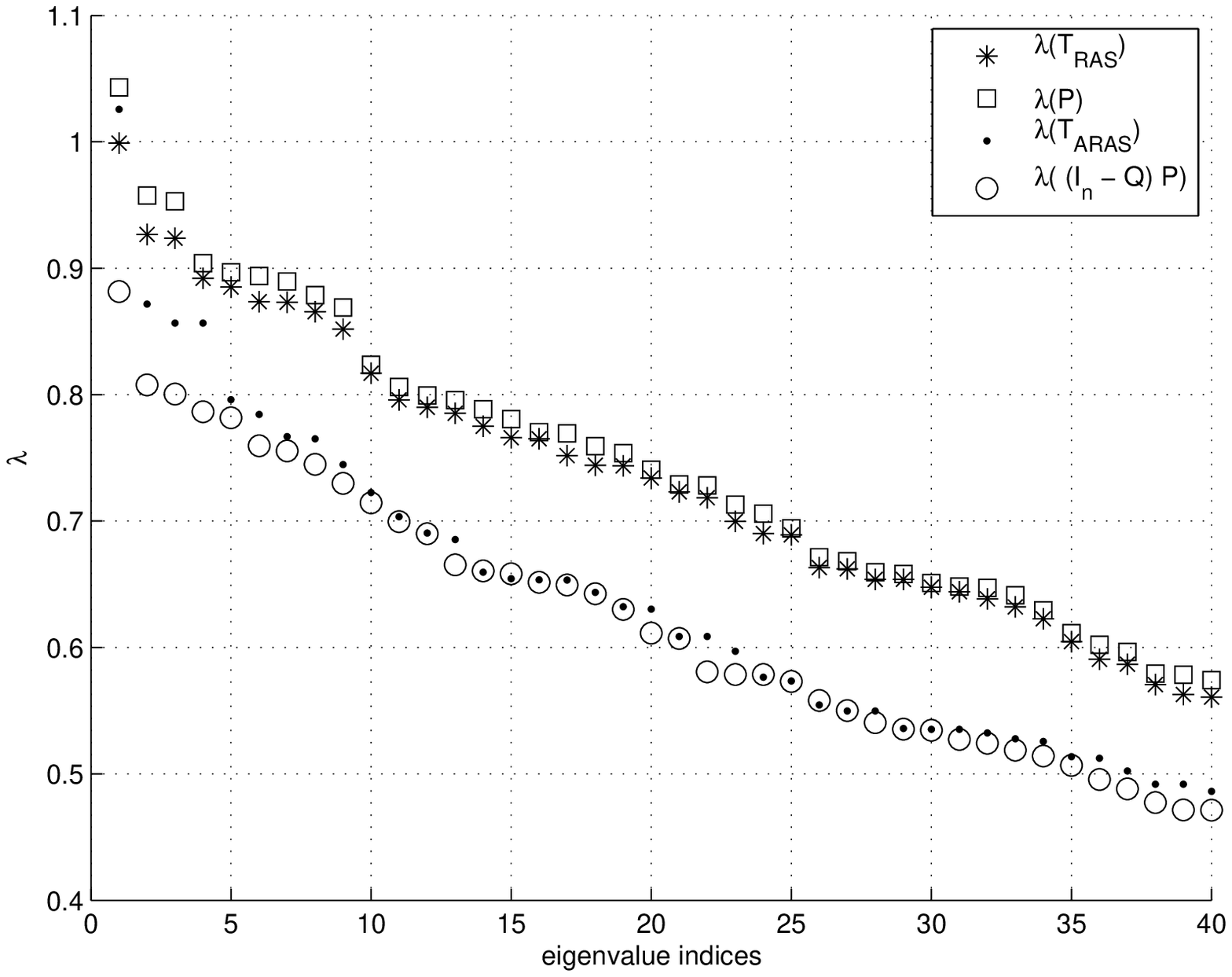}
\end{minipage}	
\caption{Eigenvalues of $T_{ARAS}(q=24)$ compared to eigenvalues of $P_{U_q}$ for a $164 \times 164$ Cartesian grid,  $p =8$,(top) Band partitioning, ARAS is built with a SVD base computed with $24$ singular vectors chosen over $2296$, (bottom) METIS partitioning, ARAS is built with a SVD base computed with $24$ singular vectors chosen over $1295$.}
\label{figure-compare-vp}
\end{center}
\end{figure}

Figure \ref{figure-compare-vp} focuses on the eigenvalue of the error transfer operators when the base is computed from SVD and both partitioning. We compute all the singular values corresponding to the number of interface points and select $24$ singular values from this set of $n$ values. For $8$ partitions with a manual partitioning, $n=2296$ and with a METIS partitioning we obtain $1295$ interfaces points.  We saw that the Aitken-RAS technique used as a Richardson iterative process diverges for the METIS partitioning. We study the spectrum of a preconditioner in the two cases and compare it to the spectrum of the RAS preconditioning method. For convenience we consider only a set of the $40$ largest eigenvalues.

The predicted values of the spectrum of the error transfer operator on the interface gives a good approximation of the spectrum of the iterative process in the case of the band partitioning. Otherwise, for the METIS partitioning, it appears that the spectrum gives a good idea of the spectrum of the iterative method but differs for the first eigenvalues. Then it appears that the first eigenvalue of $T_{ARAS}$ is greather than $1$ instead of the value computed for $\bar{Q}P$, which is less than $1$. It should explain why the iterative process diverges.

This empirical analysis shows the influence of a partitioning technique on the Schwarz preconditioner. It is important to know that if the user has the entire knowledge of the linear system he solves then he should provide a physical partitioning which can have smooth boundaries. Otherwise, if the sub-problems are  non-singular then the Schwarz method used as a preconditioning technique is efficient but can present a lack of speed in the convergence. The second point is the choice of a good base. The two bases are efficient, but the one arising from the SVD presents the best choice for time computing considerations. \\

\section{Results on industrial problems}
\label{section-results-ind}

In section \ref{section-results-academic} we pointed out that a graph partitioning to define the domain decomposition and the algorithm of approximation using the SVD of the Schwarz solutions should be a good choice to algebraically build the ARAS preconditioner. We first propose an estimation of the computing cost and then apply the preconditioner to an industrial problem.

	\subsection{Computing cost modelling}
	\label{section-cost-modelling}
	
	We want to evaluate the cost of building and applying an ARAS type preconditioner in terms of arithmetic complexity. We denote by $\cal{AC}(*)$ the arithmetic complexity of an operation $*$. 
Let considers a matrix $A \in \mathbb{R}^{m \times m}$. This matrix is decomposed into $p$ sub-matrices $A_i \in \mathbb{R}^{m_i \times m_i}$. The decomposition leads to have an interface $\Gamma$ of size $n$. The coarse interface is of size $q$. We denote by $x_{\alpha} \in \mathbb{R}^{\alpha}$ where $\alpha \in \mathbb{N}$ should be any of $m$, $m_i$, $n$ or $q$.

		\subsubsection{Arithmetic complexity of applying an ARAS type preconditioner}
Let considers the operation:
\begin{equation*}
M_{RAS}^{-1}x_m = \sum_{i=1}^{p} \tilde{R}^{T}_i A_{i}^{-1}x_{m_i} = y
\end{equation*} 

The cost of such an operation mostly consists in the $p$ operations $A_{i}^{-1}x_{m_i}$ which depends on the cost of a chosen method to inverse $A_i$ such as a Krylov methods or a LU factorization or maybe an incomplete LU factorization.

Then the cost of applying a RAS preconditioner is written as:
\begin{equation}\label{eq:RAS-apply-cost}
{\cal{AC}} \left( M_{RAS}^{-1} x_{m} \right) = p \times {\cal{AC}} \left( A_{i}^{-1} x_{m_i} \right)
\end{equation}

We now consider the operation:
\begin{equation*}
M^{-1}_{ARAS(q)} x_m = \left(I_{m}+R_{\Gamma}^{T}\mathbb{U}_q\left( \left(I_{q}-{P_{\mathbb{U}_q}}\right)^{-1}-I_{q}\right)\mathbb{U}_q^T R_{\Gamma}\right)\displaystyle\sum_{i=1}^{p}\tilde{R}_{i}^{T}A_{i}^{-1}R_{i}  x_m
\end{equation*}

The cost of such an operation consists in one application of a RAS preconditioner, the base transfer operated by $\mathbb{U_q}$ and solving $\left(I_q - P_{\mathbb{U}_q}\right) y_q = x_q$. Others or summing operations: 1 sum between 2 vectors of size $n$ and, one subtract between two vector of size $q$.

Then the cost of applying a RAS preconditioner is written as:
\begin{eqnarray*}
{\cal{AC}} \left( M_{ARAS}^{-1} x_{m} \right) &=& p \times {\cal{AC}} \left( A_{i}^{-1} x_{m_i} \right) + {\cal{AC}} \left(\mathbb{U}_q^T x_m \right) \\
& & + {\cal{AC}} \left(  \left(I_q - P_{\mathbb{U}_q}\right)^{-1} x_q \right) +  {\cal{AC}} \left(\mathbb{U}_q x_q \right) \\
& & + {\cal{AC}} \left(x_q - x_q \right)  + {\cal{AC}} \left(x_n + x_n \right) 
\end{eqnarray*}

We note that:
\begin{itemize}
\item An addition between 2 vectors of  size $n$ consists in $n$ operations. 
\item A subtraction between 2 vectors of size $n$ consists in $n$ operations. 
\item A scalar product between 2 vectors of size $n$ consists in $n$ product and $n$ sum.
\item A multiplication between a matrix with $n$ lines and $m$ columns and a vector with $m$ lines consists in $n$ scalar products of vectors of size $m$.
\end{itemize}

Then we write,
\begin{eqnarray*}
{\cal{AC}} \left(x_n + x_n \right) &=& n \\
{\cal{AC}} \left(x_q - x_q \right)  &=& q \\
{\cal{AC}} \left(\mathbb{U}_q x_q \right)  &=& m \times 2 \times q \\
{\cal{AC}} \left(\mathbb{U}^{T}_q x_m \right)  &=& q \times 2 \times m \\
\end{eqnarray*}

And,
\begin{eqnarray}
{\cal{AC}} \left( M_{ARAS}^{-1} x_{m} \right) &=& p \times {\cal{AC}} \left( A_{i}^{-1} x_{m_i} \right)  \\
& & + {\cal{AC}} \left(  \left(I_q - P_{\mathbb{U}_q}\right)^{-1} x_q \right)\\
& & + 4 \times q \times m + n + q 
\end{eqnarray}

\begin{remark}
In most cases the coarsening is such that $q \ll m$. Then the cost $${\cal{AC}} \left(  \left(I_q - P_{\mathbb{U}_q}\right)^{-1} x_q \right)$$ should be very small compared to the cost $ p \times {\cal{AC}} \left( A_{i}^{-1} x_{m_i} \right) $.
If it is the case 
\begin{equation}
{\cal{AC}} \left( M_{ARAS(q)}^{-1} x_{m} \right) = {\cal{AC}} \left( M_{RAS}^{-1} x_{m} \right) + {\cal{O}}(m)
\end{equation}

This means that the cost of one application of the ARAS preconditioner is close to the cost of one application of the RAS preconditioner when $q \ll m$.
\end{remark}

We now estimate the cost of applying an ARAS2 preconditioner:
\begin{equation}
M^{-1}_{ARAS2(q)} x_m = 2M_{ARAS(q)}^{-1} - M_{ARAS(q)}^{-1}AM_{ARAS(q)}^{-1} x_m
\end{equation}

It consists in 2 applications of ARAS(q) and one matrix vector product on the entire domain. We note that the matrix $A$ is sparse. So, denoting by $nnz(A)$ the number of non zeros elements of $A$ we write,
\begin{equation*}
{\cal{AC}}(A x_m) = 2 \times nnz(A)
\end{equation*}
Hence,
\begin{equation}
{\cal{AC}} \left( M_{ARAS2(q)}^{-1} x_{m} \right) = 2 \times {\cal{AC}} \left( M_{ARAS(q)}^{-1} x_{m} \right) + 2 \times nnz(A) + 2 {\cal{O}}(m)
\end{equation}

\begin{remark}
For $q \ll m$,
\begin{equation*}
{\cal{AC}} \left( M_{ARAS2(q)}^{-1} x_{m} \right) =2 \times {\cal{AC}} \left( M_{RAS}^{-1} x_{m} \right) + {\cal{O}}(nnz(A)) + {\cal{O}}(m)
\end{equation*}
\end{remark}

		\subsubsection{Arithmetic complexity of building the coarse space $\mathbb{U}_q$ and $P_{\mathbb{U}_q}$}

We focus here on the cost to build a base arising from the SVD of the Schwarz solutions on the interface. We refer to Algorithm \ref{algorithm-SVD-without-inverting} which proposes a robust way to implement the Aitken's acceleration without inversion.

We compute $q+2$ iterations of a RAS iterative process. It consists in applying the preconditioner on a vector $x_m$ and summing the result with another vector of the same size. We can write
\begin{equation}
{\cal{AC}}(T_{RAS} x_m) = {\cal{AC}}(M_{RAS}^{-1} x_m) + {\cal{O}}(m)
\end{equation}

Then we perform a SVD over a set of $q+2$ vectors of size $n$, $X_{q+2} \in \mathbb{R}^{n \times (q+2)}$. 
\begin{equation}
{\cal{AC}}(\text{building } \mathbb{U}_q)  \leq (q+2) \times {\cal{AC}}(T_{RAS} x_m)  + {\cal{AC}}(SVD(X_{q+2}))
\end{equation}

After this, we apply one iteration of the Schwarz iterative process on at most the $q$ first left singular vectors to build $P_{\mathbb{U_q}}$.
Thus,
\begin{eqnarray}
{\cal{AC}}(\text{building } \mathbb{U}_q  \text{ and } P_{\mathbb{U}_q})  & \leq (q+2) & \times {\cal{AC}}(T_{RAS} x_m)  \nonumber\\
& & + {\cal{AC}}(SVD(X_{q+2}) + q \times {\cal{AC}}(M_{RAS}^{-1}x_m)
\end{eqnarray}

Hence,
\begin{equation}
{\cal{AC}}(\text{building } \mathbb{U}_q \text{ and } P_{\mathbb{U}_q})  \leq 2(q+1) \times {\cal{AC}}(M_{RAS}^{-1} x_m)  + {\cal{AC}}(SVD(X_{q+2})) + {\cal{O}}(m)
\end{equation}

The cost of building the coarse space and the error transfer operator depends on the number $q$ of vectors needed. Then the ARAS(q) preconditioner will be a good choice compared to RAS if $q$ is sufficiently small compared to the number of application of the preconditioner involved in the Krylov iterative method.

\begin{remark}
This computation, following the robust algorithm \ref{algorithm-SVD-without-inverting}, is nearly two times costly than Algorithm \ref{algorithm-SVD-with-inverting} with inversion.
In fact, the building of $P_{\mathbb{U}_q}$ with Algorithm  \ref{algorithm-SVD-with-inverting} consists in inverting an error matrix of size $q$ and multiplying it, on the left, by another matrix of size $q$. For simplicity we consider those operations of order ${\cal{O}}(m)$.
\begin{equation}
{\cal{AC}}(\text{building } \mathbb{U}_q \text{ and } P_{\mathbb{U}_q})  \leq (q+2) \times {\cal{AC}}(M_{RAS}^{-1} x_m)  + {\cal{AC}}(SVD(X_{q+2})) + {\cal{O}}(m)
\end{equation}
In order to save computing, Algorithm  \ref{algorithm-SVD-with-inverting} with inversion is the best choice.
\end{remark}

		\subsubsection{Parallelization}
It is important to note that the Restricted Additive Schwarz process is fully parallel, in the sense that the inverse of $ A_i $ can be computed independently by every single process $i$ handling a domain $i$. Then for $p$ processes, we can re-write the formula \eqref{eq:RAS-apply-cost}:
\begin{equation}\label{eq:parallel-RAS-apply-cost}
{\cal{AC}} \left( M_{RAS}^{-1} x_{m} \right) ={\cal{AC}} \left( A_{i}^{-1} x_{m_i} \right)
\end{equation}

The parallelism leads to have a reduction of the matrix vector product. Then we write, for $p$ processes:
\begin{eqnarray*}
{\cal{AC}} \left(\mathbb{U}_q x_q \right)  &=& m_i \times 2 \times q \\
{\cal{AC}} \left(\mathbb{U}^{T}_q x_m \right)  &=& q \times 2 \times m_i \\
{\cal{AC}}(A x_m) & = & 2 \times nnz_{i}(A)
\end{eqnarray*}

Table \ref{table-complexity} shows the parallel arithmetic complexity in ${\cal{O}}(m_i)$ for building and applying a parallel ARAS(q) preconditioner.
\begin{table}
\begin{center}
\begin{tabular}{|c|c|c|}
\hline
Task & ${\cal{AC}}(SVD(X_{q+2}))$ & ${\cal{AC}} \left( A_{i}^{-1} x_{m_i} \right)$ \\
\hline
Building & $1$ & $q+2$ \\
\hline
Apply  & 0 & 1 \\
\hline
\end{tabular}
\caption{Parallel complexity for building and applying an ARAS(q) preconditioner}
\label{table-complexity}
\end{center}
\end{table}

\subsection{Application on a $3D$ CFD industrial case}

We consider CASE\_017 RM07 available in the sparse matrix collection \cite{UFmatrix}, which represents a $3D$ viscous case with a "frozen" turbulence. Here, the geometry is a jet engine compressor. The problem is discretized among 54527 nodes. Seven variables per node are considered. The resulting matrix is of size $381 689$ with $37 464 962$ non-zeros. The matrix is not symmetric.

We use a PETSc-MPI implementation with a PARMETIS partitioning taking into account the block structure and the weight of each block. The matrix is partitioned in four parts, with a minimum overlap. We apply the ARAS and ARAS2 left-preconditioners with 40 basis vectors computed following Algorithm \ref{algorithm-SVD-without-inverting}.  

\begin{figure}[t]
\begin{center}
\begin{minipage}{13cm}
	\includegraphics[height=80mm,width=60mm]{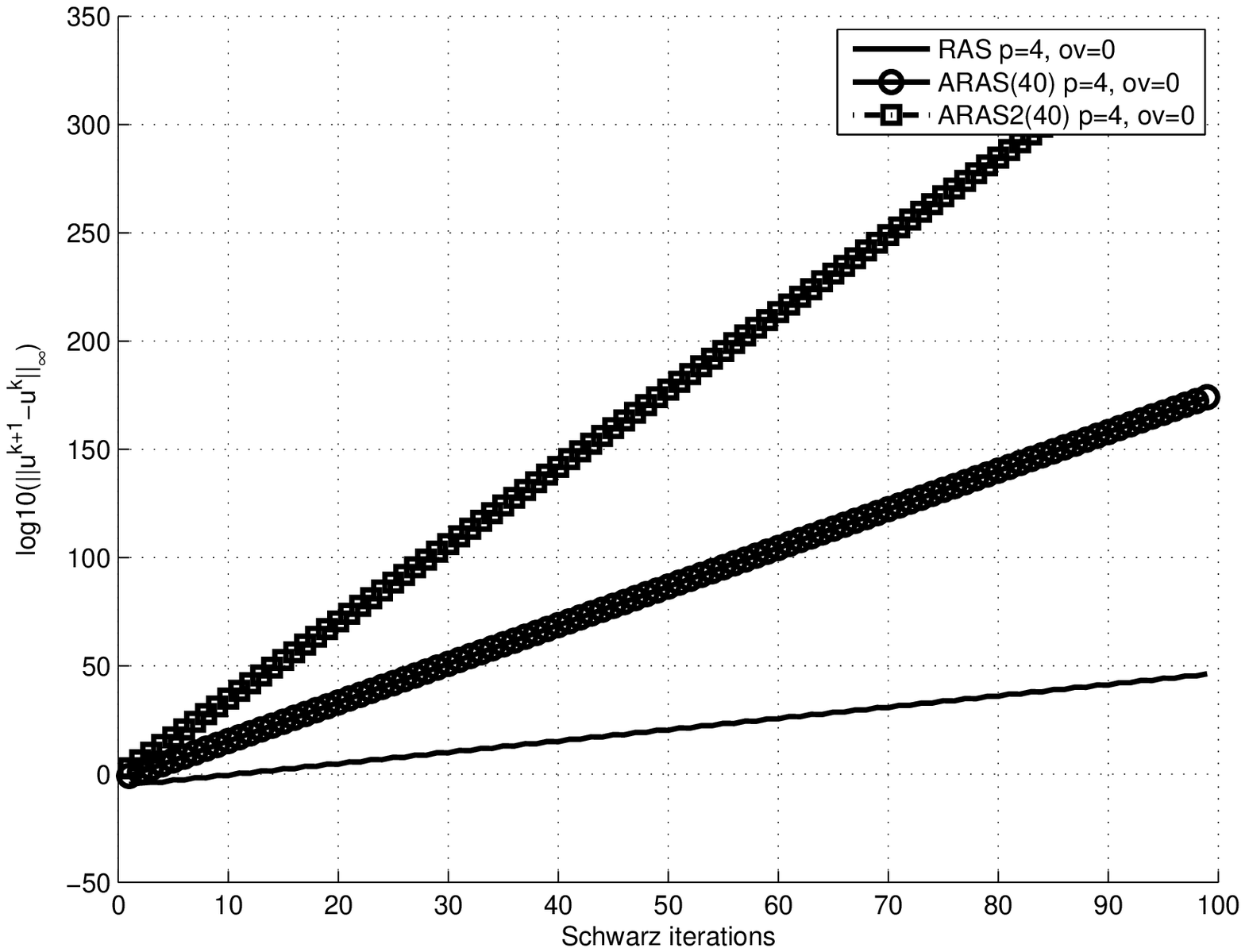}
	\includegraphics[height=80mm,width=60mm]{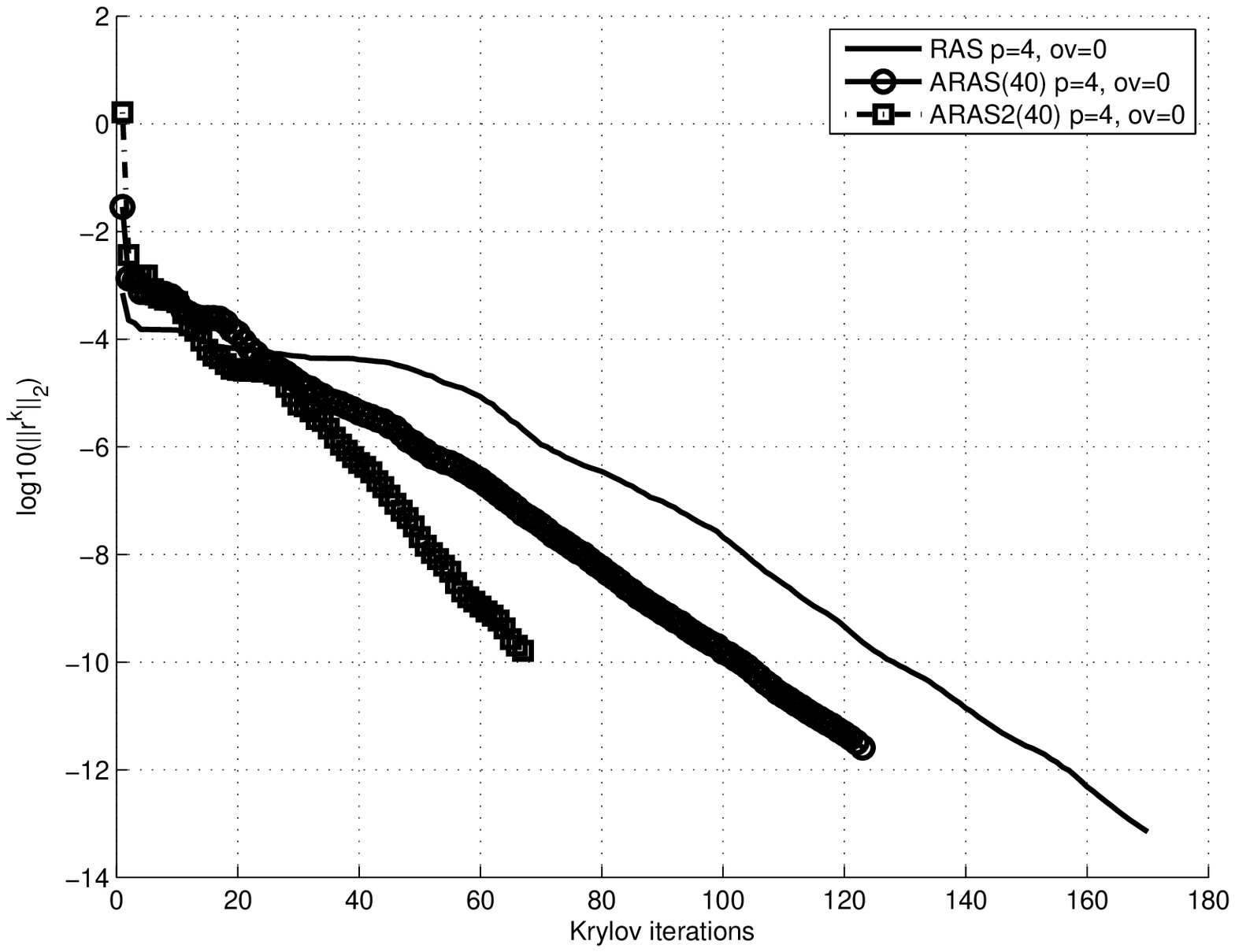}
\end{minipage}	
\caption{Solving 3D Navier Stokes equation  (CASE RM07), PARMETIS partitioning with weights, $p =4$, overlap $0$, ARAS2 is built with a SVD basis, (left) Convergence of Iterative  Schwarz Process, (right) convergence of GMRES method preconditioned by RAS and ARAS (built with $tol = 10^{-6}$).}
\label{figure-3D-results}
\end{center}
\end{figure}

Figure \ref{figure-3D-results} shows the convergence of the preconditioned GMRES and the convergence of the Richardson processes associated to the preconditioners. In order to disminish the cost of building phase, the tolerance is set to $10^{-6}$ for the local solution, but only for the building phase.Then, the iterations to build $P_{\mathbb{U}_q}$ take less cpu-time and memory allocation than one application of the preconditioner during the solution phase. The ARAS(q) preconditioner has been presented as a left preconditioner. The stopping criteria used for the GMRES method implemented in PETSc is based on the relative residual. This is why the curves are not going to the same tolerance. We discuss this point with the following results. For a minimal overlap, without knowledge of the underlying equations and mesh design, the ARAS preconditioner is efficient, and also its mutliplicative version, ARAS2. \\

Table \ref{table-3D-results} shows the performance results corresponding to Figure \ref{figure-3D-results} on an SGI Altix Xe340.  While applying ARAS(40) or ARAS2(40), both Solution Time and Memory allocation are reduced. The time involved in the building phase of the preconditioner depends of the choice of local factorization and the kind of algorithm chosen. The different times follows the formula of the arithmetic complexity presented in Subsection \ref{section-cost-modelling}.
The relative residual for RAS and ARAS(q) is the same. We explain the difference of convergence between ARAS and ARAS2 for this case by the fact that the Richardson processes diverge and leads to an amplified inaccurary in ARAS2.

\begin{table}[h]
\begin{center}
\begin{tabular}{|c|c|c|c|c|}
\hline
 Prec. & Building & Solution & Max. Loc. Mem & $||Ax - b||_2/||b||_2$ \\
           & Time (s.)     & Time (s.)     & (M.O.)      &                                      \\
\hline
 RAS   & 8.684           &  1552.89 & 1068          & 1.5704 e-09                \\
 ARAS (40)  & 429.943  &  1086.63 &  1048   & 9.7492e-10                  \\
 ARAS2 (40)  & 446.454  &  1174.97 &  1010   &   2.2760e-07   \\
 \hline
\end{tabular}
\caption{Solving 3D Navier Stokes equation  (CASE RM07), PARMETIS partitioning with weights, $p =4$, overlap $0$, ARAS2 is built with a SVD basis, GMRES method preconditioned by RAS and ARAS (built with $tol = 10^{-6}$).}
\label{table-3D-results}
\end{center}
\end{table}

\section{Conclusion}
\label{section-outlook}
We presented an integration of the Aitken acceleration technique in the RAS preconditioning. This integration leads to a multi-level preconditioner. One level is the entire domain, while the second is the entire artificial interface. Since the computation of the error transfer operator is costly, we propose an algebraic computation of a coarse space, built from the SVD decomposition of Schwarz solutions on the interface. The results is a cheap and fully algebraic enhancement of the RAS preconditioner. 
An analysis of the convergence of this preconditioner is given when the basis is built analytically and shows the effect of
 the preconditioner depending on the choice of the mode to be accelerated.
 Finally, a result is provided on a $3D$ industrial case without knowledge of the underlying equations and the mesh design.\\
Future work should focus on other algorithms to build algebraically the Aitken acceleration in order to reduce the time spent in the building time. Another issue concerns work on partitioning techniques for domain decomposition. Actually, there is no technique to provide local system with insurance of inversion.
 
\section*{Acknowledgements}
This work was funded by the French National Agency of Research under the contract ANR-TLOG07-011-03 LIBRAERO. 
The work of the second authors was also supported by the r\'egion Rh\^one-Alpes through the project CHPID of 
the cluster ISLE.\\
Authors are grateful to FLUOREM for providing the industrial test cases and the PETSc code setting the problem solving environment.\\
The experiments were done on the cluster SGI-XEON of the Centre pour le D\'eveloppement du Calcul Scientifique Parall\`ele of the Universit\'e Lyon 1.

\topsep=0.0ex
\parsep=0.0ex
\parskip=0.0ex
\itemsep=0.0ex

\bibliographystyle{siam} 
\def\cprime{$'$}

\end{document}

%% file: figures/figure-2D-domain-poisson.tex
\begin{tikzpicture}
\begin{scope}
	\node at (0,0) [name=A]{};
	\node at (0,5) [name=B]{};
	\node at (12,5)[name=C]{};
	\node at (12,0)[name=D]{};
	
	\draw (A.center)--(B.center)--(C.center)--(D.center)--(A.center);

	\node at (5.5,5) [name=E]{};
	\node at (6.5,5) [name=F]{};	
	\node at (6.5,0) [name=G]{};
	\node at (5.5,0) [name=H]{};		

	\draw (E) node[above] {$\Gamma_2$};
	\draw (F) node[above] {$\Gamma_1$};
	\draw [-,very thick] (E.center)--(H.center);
	\draw [-,very thick] (F.center)--(G.center);
	
	\node at (3,2.5) [name=O1]{$\Omega_1$};
	\node at (9,2.5) [name=O2]{$\Omega_2$};
	
	\draw (0,5/pi) node[left]{$y$};
	\draw (12,0) node[below]{$x$};
	\draw [->,thick] (A.center)--(0,5/pi);
	\draw [->,thick] (A.center)--(12,0);
\end{scope}
\end{tikzpicture}